\documentclass{article}
\usepackage{inputenc,amsfonts,amsmath,amssymb,mathtools}
\usepackage[all]{xy}
\usepackage{amsthm}
\usepackage{hyperref}
\usepackage[left=2.5cm, right=2.5cm, top=2cm]{geometry}
\usepackage{extarrows}
\usepackage{graphicx,subcaption}
\usepackage{wrapfig}
\usepackage{multirow}

\usepackage{aliascnt}
\newcommand{\thdef}[2]{
	\newaliascnt{#1}{theorem}  
	\newtheorem{#1}[#1]{#2}
	\aliascntresetthe{#1}  
	\newtheorem*{#1*}{#2}
	\expandafter\newcommand\expandafter{\csname #1autorefname\endcsname}{#2}
}

\newcommand{\edge}{\,
	\begin{tikzpicture}
		\draw (0,.5ex) -- (2ex,.5ex);
		\draw[fill] (0,.5ex) circle (.2ex);
		\draw[fill] (2ex,.5ex) circle (.2ex);
		\draw[fill] (0,0) circle (0);
	\end{tikzpicture}\,
}

\usepackage{tikz-cd,pgfplots}
\usetikzlibrary{angles,quotes}
\title{Elliptic log symplectic brackets on projective bundles}

\author{Mykola Matviichuk\footnote{Imperial College London, m.matviichuk@imperial.ac.uk} }
\date{}

\DeclareMathOperator*{\res}{res}

\newtheorem{theorem}{Theorem}[section]
\thdef{thm}{Theorem}
\thdef{lemma}{Lemma}
\thdef{proposition}{Proposition}
\thdef{corollary}{Corollary}
\thdef{conjecture}{Conjecture}
\thdef{problem}{Problem}
\theoremstyle{definition}
\thdef{definition}{Definition}
\thdef{example}{Example}
\thdef{remark}{Remark}

\thdef{theoremAlph}{Theorem}

\thdef{definitionAlph}{Definition}

\thdef{conjectureAlph}{Conjecture}

\usepackage{todonotes}

\newcommand{\X}{\mathsf{X}}
\newcommand{\D}{\mathsf{D}}
\newcommand{\Pf}{{\rm Pf}}
\newcommand{\ii}{{\rm i}}
\newcommand{\n}{n}

\newcommand{\sageNgon}[3][0.4]{
\pgfdeclarelayer{background layer}
\pgfdeclarelayer{foreground layer}
\pgfsetlayers{background layer,main,foreground layer}
\begin{tikzpicture}[baseline=-1ex,
rotate=360/2/#2,
scale=#1,line join = round 
    ] 

\begin{pgfonlayer}{main}
\foreach \n in {1,...,#2}
{
    \coordinate (v\n) at ({-90+((\n-1)*360/#2)}:1.5);
}

\foreach \n in {1,...,#2}
{
	\foreach \m in {1,...,#2}
	\draw (v\n) -- (v\m);
}
\end{pgfonlayer}

#3
\begin{pgfonlayer}{foreground layer}
\foreach \n in {1,...,#2}
{
	\draw[fill] (v\n) circle (0.03/#1);
}
\end{pgfonlayer}
 \end{tikzpicture}
}

\newcommand{\smoothedge}[5][blue]{
\begin{scope}
\def\Ea{#2}
\def\Eb{#3}
\def\Ra{#4}
\def\Rb{#5}
\draw[very thick,#1] (\Ea) -- (\Eb);
\foreach \x in {\Ra,\Rb} {
		\begin{pgfonlayer}{background layer}
		\begin{scope}
			\clip (\Ea) -- (\x) -- (\Eb);
			
			\draw [draw,opacity=0.5,fill,#1] (\x) circle (0.33) ;
		\end{scope}
		\end{pgfonlayer}
}
\end{scope}
}

\begin{document}

\maketitle

\begin{abstract}
Let $\X$ be the product of a complex projective space  and a polydisc. We study Poisson brackets on $\X$ that are log symplectic, that is, generically symplectic and such that the inverse $2$-form has only first order poles. We propose a method of constructing such Poisson brackets that additionally are elliptic, in a precise sense. Our method relies on the local Torelli theorem for log symplectic manifolds of Pym, Schedler and the author, and uses combinatorics of smoothing diagrams. We demonstrate effectiveness of the method on a series of examples, recovering, in particular, all log symplectic cases of elliptic Feigin-Odesskii Poisson brackets $q_{n,k}$ on $\mathbb{P}^{n-1}$.
\end{abstract}

\section*{Introduction}

Elliptic Poisson algebras, together with their quantizations, have earned significant attention from both mathematicians and physicists, as they are considered among the most intricate examples of Poisson algebras and exhibit remarkably rich geometric properties. Elliptic algebras of Feigin and Odesskii \cite{Odesskii2002} are directly linked to some of the most complex solutions of the Yang-Baxter equation \cite{Chirvasitu2023}. Elliptic zastava spaces \cite{Finkelberg2023} are expected to be linked to a generalization of the BFN construction \cite{Finkelberg2018} of the Coulomb branch of a quiver gauge theory \cite[Section 4]{Finkelberg2023}.

The goal of the current article is to propose a systematic way of constructing new examples of elliptic Poisson algebras, using the deformation theory developed in \cite{Matviichuk2020}. Everywhere below, we will work in the realm of log symplectic Poisson manifolds, since this is the setup of \cite{Matviichuk2020}. The key idea is to start with a log symplectic manifold of a simple kind, namely, such that the polar divisor of the log symplectic form has only normal crossings singularities, and deform it to obtain a more complicated log symplectic manifold. 

One of the easiest and, perhaps, most natural cases where such an approach applies is when the complex manifold carrying the log symplectic form is $\X = \mathbb{P}^{\n-1}$, the complex projective space of even dimension $\n-1$. An obvious source of log symplectic forms on $\X$ are the toric ones, which are invariant under the action of the torus $(\mathbb{C}^*)^{\n-1}$. These have simple poles precisely along each coordinate hyperplane. The periods of a toric log symplectic form $\omega_0$ on $\X$ can be conveniently organized in an $n\times n$ skew-symmetric matrix of complex numbers $B=(b_{ij})_{i,j=0}^{n-1}$ whose rank is $n-1$ and whose rows sum to zero. Following \cite{Matviichuk2020}, we incorporate the periods of $\omega_0$ into a certain graph with decorations, called \textit{smoothing diagram}. 
The graph has $\n$ vertices, one for each irreducible component of the polar divisor. Every two vertices $i$, $j$ in the graph are connected by an edge, reflecting the fact that every two irreducible components have non-trivial intersection.
Furthermore, we color an edge $i\edge j$ and call it \textit{smoothable} if $b_{ij}\not=0$ and for each $k\not=i,j$ the quotient
$$
\frac{b_{jk}+b_{ki}}{b_{ij}} \text{ is a non-negative integer.}
$$
The presence of a smoothable edge guarantees (\autoref{prop:smoothable_edge_deform} for  $m=0$) that there exists a Poisson deformation of $\pi_0=\omega_0^{-1}$ that smoothes out the corresponding irreducible component of the singular locus of the polar divisor, hence the terminology. Remarkably, the Poisson deformations given by different smoothable edges are jointly unobstructed (\autoref{thm:deforming_several_smootables} for  $m=0$).

The results above prompt the combinatorial question of identifying the potential types of subgraphs that the smoothable edges can form within the complete graph on $n$ vertices. One can check (\autoref{lm:pos_rat_valency}) that in fact each connected component of such a subgraph has to be either a chain, or a cycle. Our working hypothesis is that the presence of cycles of smoothable edges indicates that the geometry of the corresponding Poisson deformation has something to do with an ellitptic curve. We will make this more precise in a moment.

While analyzing the combinatorics of the cycles of smoothable edges, as we do in \autoref{sec:smoothableCycles}, one can notice that the condition on the rank of the matrix $B$ needs an extra care and that it is easier to carry out the classification of such cycles if the condition on the rank is dropped. This observation led us to considering a slightly more general complex manifold $\X = \mathbb{P}^{n-1}\times \mathbb{D}^m$,  where $\mathbb{D}^m\subset \mathbb{C}^m$ is a polydisc. This manifold admits log symplectic forms that realize all possible skew-symmetric matrices $B$ whose rows sum to zero, as long as $m$ is chosen large enough (see \autoref{lm:boundDimPolydisc}). Notions of smoothing diagram and smoothable edges generalize verbatim (see \autoref{sec:deformTheory}) to log symplectic forms on $\X = \mathbb{P}^{n-1}\times \mathbb{D}^m$ that are semi-toric, i.e. invariant under the infinitesimal action of $(\mathbb{C}^*)^{n-1}\times \mathbb{C}^m$.

We remark that such an $\X$ may be considered as a projective bundle over $\mathbb{D}^m$. By splitting off symplectic direct summands, one can arrange the fibers of this bundle to be coisotropic. This setup allows us to define a Higgs field on the deprojectivized vector bundle capturing the mixed terms of the Poisson tensor, a construction that was studied in detail in \cite{MatviichukPhD}. In the examples below, the Higgs field remains unchanged under the Poisson deformation and plays a crucial role in the analysis of the deformed log symplectic bracket.

By an elliptic curve we mean a compact Riemann surface of genus one. We propose the following

\begin{definitionAlph}\label{def:elliptic}
A log symplectic manifold is called \textit{elliptic} if the space of its codimension two symplectic leaves is an elliptic curve.
\end{definitionAlph}

\begin{conjectureAlph}\label{conj:ellipticVSsmoothableCycle}
Let $\omega_0$ be a semi-toric log symplectic form on $\X=\mathbb{P}^{n-1}\times \mathbb{D}^m$, $n\ge 3$, $m\ge 0$. Let $\pi$ be a Poisson deformation of $\pi_0=\omega_0^{-1}$ produced, in the sense of \autoref{thm:deforming_several_smootables}, from a  cycle of $n$ smoothable edges. Then $\pi$ is elliptic.

Conversely, every elliptic log symplectic Poisson bracket on $\X=\mathbb{P}^{n-1}\times \mathbb{D}^m$ can be obtained this way.
\end{conjectureAlph}


We provide evidence for the conjecture in \autoref{sec:deformTheory} (the running example), \autoref{subsec:C41deform}, \autoref{subsec:X5deform}, where we construct seemingly new examples of elliptic Poisson brackets on $\mathbb{P}^3\times \mathbb{D}^1$ and $\mathbb{P}^4\times\mathbb{D}^2$, and in \autoref{subsec:FO}, where we demonstrate that the log symplectic Feigin-Odesskii brackets $q_{n,k}$ on $\mathbb{P}^{n-1}$ can be obtained from cycles of smoothable edges.

Assuming the conjecture is true, one can target classification of elliptic log symplectic brackets on $\mathbb{P}^{n-1}\times \mathbb{D}^m$ using the following combinatorial result.

\begin{theoremAlph}\label{thm:smoothableCycles}
Let $B$ be an $n\times n$ skew-symmetric matrix whose rows sum to zero. Then $B$ has a cycle of $n$ smoothable edges if and only if, perhaps after rescaling and conjugation by a permutation matrix, it appears in the following list of matrices (all defined in \autoref{sec:smoothableCycles})
\begin{itemize}
\item $\mathcal{C}_{\n,k,I}$, where $1\le k< \n/2$, and $I$ is a sequence of zeros and ones of length $d=\gcd(\n,k)$,
\item $\mathcal{X}_4$, $\mathcal{X}_5$,
\item $\mathcal{Y}_{2(2k+1)}$, $k\ge 1$,
\item $\mathcal{Z}_{5(2k+1)}$, $k\ge 1$.
\end{itemize}

\end{theoremAlph}



We verify that the period matrices in this list that can be realized on $\mathbb{P}^{n-1}$ are precisely the ones that lead to log symplectic Feigin-Odesskii  brackets $q_{n,k}$ (\autoref{thm:FOviaSmoothableCycles}). This suggests that perhaps all elliptic log symplectic brackets on $\mathbb{P}^{n-1}$ are Feigin-Odesskii.

It is expected that some of the elliptic log symplectic brackets on $\mathbb{P}^{n-1}\times \mathbb{D}^m$ have quantizations that can be constructed via elliptic dynamical $R$-matrices \cite{Schiffmann1998,Etingof2000}.

\textbf{Acknowledgements.} I would like to thank Brent Pym and Travis Schedler for their unwavering support, guidance, and valuable insights throughout this project. I also thank Marco Gualtieri for engaging and clarifying discussions and Pavel Etingof for suggesting the intriguing connection to dynamical $R$-matrices.

\newpage

\tableofcontents

\section{Log symplectic manifolds}

By a \textit{log symplectic} manifold we mean a complex manifold $\X$ of even dimension $d$ equipped with a closed meromorphic $2$-form $\omega$ that has only a first order pole along a reduced divisor $\D\subset \X$, and such that the Pfaffian $\Pf(\omega) = \omega^{d/2}$ is a non-trivial meromorphic volume form on $\X$ having a first order pole along $\D$. If $(\X,\D,\omega)$ is log symplectic, then the bivector $\pi = \omega^{-1}$ is holomorphic Poisson, and its Pfaffian $\Pf(\pi)$ is a holomorphic co-volume form that has a first order zero along $\D$. From now on, we may refer to the divisor $\D$ as either the polar divisor of $\omega$, or the degeneracy divisor of $\pi$. Note that $\omega$ defines a class $[\omega]$ in the de Rham cohomology ${\rm H}^2_{dR}(\X\setminus \D)$. If $\D$ has only normal crossings singularities, then locally near $p\in \D$, the class $[\omega]\in {\rm H}^2_{dR}(\X \setminus \D)$ is given by the collection of complex numbers $b_{ij}$, $0\le i,j\le {\n-1}$, called the \textit{biresidues} of $\omega$ at $p$, where $\n$ is the number of local irreducible components of $\D$ at $p$. Concretely, denoting $\D_0$,...,$\D_{\n-1}$ the local irreducible components of $D$, we take the residue twice to obtain
$$
b_{ij} = \res_{\D_j\cap \D_i} \res_{\D_i} \omega
$$
for $i\not=j$, and $b_{ii} = 0$. Note that $b_{ji} = -b_{ij}$. If we choose local coordinates $x_0,...,x_{d-1}$ on $\mathcal{U}\subset \X$ centered at $p$ so that $\D_i  \cap \mathcal{U}=\{x_i=0\}$, $0\le i \le {\n-1}$, then $\omega$ will have an expression
$$
\omega = \sum_{0\le i<j\le {\n-1}} b_{ij} \dfrac{dx_i}{x_i} \wedge \dfrac{dx_j}{x_j} + \sum_{i=0}^{\n-1} \sum_{j=0}^{d-1} f_{ij} \dfrac{dx_i}{x_i} \wedge dx_j +  \sum_{0\le i<j\le= d-1} g_{ij} dx_i \wedge dx_j,
$$
where $g_{ij}$, $h_{ij}$ are locally defined holomorphic functions. A Moser argument shows that the biresidues alone are enough to recover $\omega$ locally near $p$ up to an analytic isomorphism \cite[Theorem 3.7, Example 3.9]{Pym2017}.

Let $\X=\mathbb{P}^{\n-1} \times \mathbb{D}^m$, $\n\ge2$, $m\ge 0$, where $\mathbb{P}^{\n-1}$ is the complex projective space of dimension $\n-1$ and $\mathbb{D}^m$ is a polydisc in $\mathbb{C}^m$. We assume everywhere below that  $\n-1+m$ is even. We will use the notation $\partial\, \mathbb{P}^{\n-1}$ for the toric anti-canonical divisor $\cup_{i=0}^{\n-1} {\mathsf H}_i$, where each ${\mathsf H}_i = \{[y_0\colon ... \colon y_{\n-1}]\in\mathbb{P}^{\n-1}:y_i=0\}$ is a coordinate hyperplane. We will denote by $\partial \X$ the divisor $\D = (\partial\,\mathbb{P}^{\n-1})\times \mathbb{D}^m$, whose irreducible components are $\D_i = {\mathsf H}_i \times \mathbb{D}^m$, $i=0,1,...,\n-1$.
Note that the intersection of any two irreducible components of $\partial \X$ is connected, so the biresidue $b_{ij}$ is independent of the choice of a point $p\in \D_i \cap \D_j$. In what follows, we also choose the standard coordinates $x_1,...,x_m$ on $\mathbb{D}^m\subset\mathbb{C}^m$ and homogeneous coordinates $y_0,...,y_{\n-1}$ on $\mathbb{P}^{\n-1}$. Moreover, we will often use affine coordinates on $\mathbb{P}^{\n-1}$, such as $z_k=\frac{y_k}{y_0}$, $k=1,...,\n-1$. 

\begin{lemma}
Let $\omega$ be a log symplectic form on $\X=\mathbb{P}^{\n-1}\times \mathbb{D}^m$ whose polar divisor is $\partial \X$. Then its biresidues $(b_{ij})_{i,j=0}^{\n-1}$ satisfy
\begin{equation}\label{eq:biresidues_add_up_to0}
\sum_{j=0}^{\n-1} b_{ij} = 0, ~~\text{for each}~ i.
\end{equation}
\end{lemma}

\begin{proof}
Without loss of generality, we can assume $i>0$.  Consider the closed logarithmic $1$-form 
$$\alpha_i = \res\limits_{z_i=0} \omega =
\sum_{j=1}^{\n-1} b_{ij} ~\dfrac{dz_j}{z_j} + \sum_{\ell=1}^m f_\ell ~dx_\ell,
$$
where each $f_\ell$ is a holomorphic function in $x$-variables. Changing the coordinates $z_1,...,z_{\n-1}$ to another set of affine coordinates, say $\widetilde{z}_0=\frac{1}{z_{\n-1}}$, $\widetilde{z}_j = \frac{z_j}{z_{\n-1}}$, $j=1,...,\n-2$, we obtain
$$
\alpha_i = - \sum_{j=1}^{n-1} b_{ij}~\dfrac{d\widetilde{z}_0}{\widetilde{z}_0}+ \sum_{j=1}^{\n-2} b_{ij} ~\dfrac{d\widetilde{z}_j}{\widetilde{z}_j} + \sum_{\ell=1}^m {f}_\ell ~dx_\ell = \sum_{j=0}^{\n-2} b_{ij} ~\dfrac{d\widetilde{z}_j}{\widetilde{z}_j} + \sum_{\ell=1}^m {f}_\ell ~dx_\ell,
$$
which implies that $b_{i0} = -\sum_{j=1}^{\n-1} b_{ij}$.
\end{proof}

For concrete calculations, it is convenient to consider the following type of log symplectic forms on $\X$. A log symplectic form $\omega$ on $\X$ is called \textit{semi-toric}, if it satisfies ${\rm Lie}_{z_k\partial_{z_k}}(\omega)=0$, $k=1,...,{\n-1}$, and ${\rm Lie}_{\partial_{x_j}}(\omega)=0$, $j=1,...,m$.

\begin{lemma}\label{lm:boundDimPolydisc}
Let $\n\ge3$, and $B = (b_{ij})_{i,j=0}^{\n-1}$ be a skew-symmetric matrix whose rows sum to zero. Then there exists a log symplectic form $\omega$ on $\X = \mathbb{P}^{\n-1} \times \mathbb{D}^m$ with the toric polar divisor $\partial X$ and biresidues $(b_{ij})_{i,j=0}^{\n-1}$ if and only if ${\n-1}+m$ is even and ${\rm rank}(B)\ge {\n-1}-m$. In such a case, $\omega$ can be chosen to be semi-toric.
\end{lemma}

\begin{proof}
Necessity of the inequality ${\rm rank}(B)\ge {\n-1}-m$ follows by applying \cite[Lemma 3.8]{Matviichuk2020} to a minimal stratum, say the stratum $[1:0:\dots:0] \times \mathbb{D}^m$.

Let us show sufficiency. By adding the symplectic direct summands $(\mathbb{D}^2,\omega_{\rm st}=dx_1\wedge dx_2)$, we can reduce the proof to the case when ${\rm rank}(B)={\n-1}-m$. Note that due to the property \eqref{eq:biresidues_add_up_to0}, the submatrix $(b_{ij})_{i,j=1}^{\n-1}$ has the same rank as $B$. Using skew-symmetry of $B$, we can choose a subset $J\subset \{1,2,...,\n-1\}$ of cardinality $\n-1-m$ such that the submatrix $(b_{ij})_{i,j\in J}$ also has rank $\n-1-m$. Let $0\le i_1<...<i_m\le \n-1$ be the complementary indices to $J$. Then we can define the desired semi-toric log symplectic form $\omega$  as 
$$
\omega = \sum_{\substack{i,j\in J\\i<j}} b_{ij}~ \dfrac{dz_i}{z_i} \wedge \dfrac{dz_j}{z_j} +    \sum_{k=1}^m \dfrac{dz_{i_k}}{z_{i_k}} \wedge dx_k.
$$
\end{proof}

Two log symplectic manifolds $(\mathsf{Y}_1,\omega_1)$  and $(\mathsf{Y}_2,\omega_2)$ are called \textit{stably isomorphic}, if there are symplectic manifolds $(\mathsf{S}_1,\varpi_1)$ and $(\mathsf{S}_2,\varpi_2)$ such that $(\mathsf{Y}_1\times \mathsf{S}_1,\omega_1+\varpi_1)$ is isomorphic to $(\mathsf{Y}_2\times \mathsf{S}_2,\omega_2+\varpi_2)$.

\begin{lemma}
Let $\omega_k$ be a log symplectic form on $\X_k = \mathbb{P}^{\n-1}\times \mathbb{D}^{m_k}$ with the toric polar divisor $\partial \X_k$, $k=1,2$. Let $\omega_1$ and $\omega_2$ have identical collections of biresidues. Then $(\X_1,\omega_1)$ is stably isomorphic to $(\X_2,\omega_2)$, perhaps after shrinking the polydiscs $\mathbb{D}^{m_k}$, $k=1,2$.
\end{lemma}

\begin{proof} By adding symplectic direct summands, one can arrange $m_1=m_2=m$. In this case, a Moser argument (e.g. \cite[Lemma 3.6]{Pym2017}) shows that $(\X_1,\omega_1)$ is  isomorphic to $(\X_2,\omega_2)$, perhaps after shrinking the polydiscs. The time-dependent vector field can be chosen to be tangent to the central fiber $\mathbb{P}^{\n-1} \times (0,..,0)$, so that its time one flow is well-defined in the neighborhood of this fiber.
\end{proof}

Note that every Poisson bivector on $\mathbb{P}^{\n-1}\times \mathbb{D}^m$ is the projectivization of a Poisson bivector on $\mathbb{C}^{\n}\times \mathbb{D}^m$ \cite[Theorem 3.1.5]{MatviichukPhD}. We will often use homogeneous coordinates to write down a Poisson bivector $\pi$ on $\X=\mathbb{P}^{\n-1}\times \mathbb{D}^m$, even though technically such an expression would represent a  Poisson bivector on $\mathbb{C}^{\n}\times \mathbb{D}^m$. The degeneracy divisor of $\pi$ on $\X$ in this notation is given by the zero locus of $\sum\limits_{i=0}^{n-1} y_i\partial_{y_i}\wedge \pi^{ {\dim\X}/{2}}$.

\newpage
\section{Recap on the deformation theory of log symplectic manifolds}\label{sec:deformTheory}

In \cite{Matviichuk2020}, the authors described the deformations of a log symplectic manifold $(\X,\D,\omega)$, whose polar divisor $\D$ has only normal crossings singularities. Note that the deformed polar divisor is allowed to have more complicated singularities in the setup of \cite{Matviichuk2020}. Let us discuss in detail the  case when $\X = \mathbb{P}^{\n-1} \times \mathbb{D}^m$, and $\D$ is the toric divisor $ \partial \X = (\partial\, \mathbb{P}^{\n-1})\times \mathbb{D}^m$.

\begin{figure}[h]
\centering
\begin{subfigure}[t]{0.3\textwidth}
\begin{tikzpicture}[scale = 1, baseline = -1.5]
\coordinate (a) at (1,0);
\coordinate (b) at (3,1);
\coordinate (c) at (3,0);

\clip (-0,-0.5) rectangle (4,1.25);
\begin{scope}
			\clip (a) -- (c) -- (b);
\draw [draw,opacity=1.1,fill,blue] (a) circle (0.7) ;
\end{scope}

\draw[ thick]  (a) node[left]{$k$} --  (b) node[right] { $i$};
\draw[ thick]  (a) --  (c) node[ right] { $j$};
\draw[thick, blue] (b)--(c);

\draw[fill] (a) circle (0.03);
\draw[fill] (b) circle (0.03);
\draw[fill] (c) circle (0.03);
\end{tikzpicture}
\caption{Darkly colored angle, $\theta_{ijk}=2$}
\end{subfigure}
\hskip1cm
\begin{subfigure}[t]{0.3\textwidth}
\begin{tikzpicture}[scale = 1, baseline = -1.5]
\coordinate (a) at (1,0);
\coordinate (b) at (3,1);
\coordinate (c) at (3,0);

\clip (-0,-0.5) rectangle (4,1.25);
\begin{scope}
			\clip (a) -- (c) -- (b);
\draw [draw,opacity=0.4,fill,blue] (a) circle (0.7) ;
\end{scope}

\draw[ thick]  (a) node[left]{$k$} --  (b) node[right] { $i$};
\draw[ thick]  (a) --  (c) node[ right] { $j$};
\draw[thick, blue] (b)--(c);

\draw[fill] (a) circle (0.03);
\draw[fill] (b) circle (0.03);
\draw[fill] (c) circle (0.03);
\end{tikzpicture}
\caption{Lightly colored angle, $\theta_{ijk}=1$}
\end{subfigure}

\caption{Two types of colored angles}
\label{fig:colored_angles_types}
\end{figure}
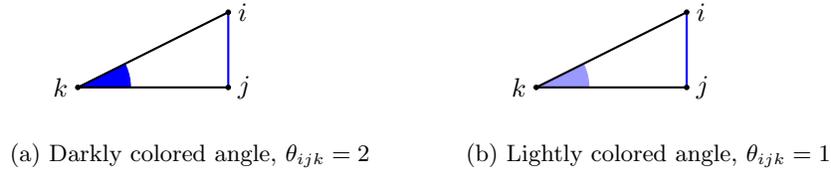

To a log symplectic form $\omega_0$ on $\X$ with polar divisor $\D=\partial \X$, we associate a combinatorial gadget called \textit{smoothing diagram}, defined as follows. First, we construct a graph, which has $\n$ vertices, one for each irreducible component of $\D$. Every two vertices $i$, $j$ in the graph are connected by an edge, signifying the fact that every two irreducible components of $\D$ have non-trivial intersection. Certain edges and angles are colored according to the following recipe that involves the biresidues $b_{ij}$, $0\le i,j\le {\n-1}$, of the log symplectic form $\omega_0$. An edge $i\edge j$ is called \textit{smoothable}, if $b_{ij}\not=0$, and the number
$$
\theta_{ijk} := \dfrac{b_{jk}+b_{ki}}{b_{ij}}
$$
is a non-negative integer, whenever $k\not=i,j$. In the smoothing diagram, we color each smoothable edge. Additionally, we color each angle $i\edge k\edge j$
%
%
%
such that $i\edge j$ is smoothable and $\theta_{ijk}$ is non-zero. Note that the equations \eqref{eq:biresidues_add_up_to0} imply that, if $b_{ij}\not=0$, then $\sum_{k\not=i,j} \theta_{ijk} = 2$. Therefore, for a smoothable edge $i\edge j$, most indices $k\not= i,j$ satisfy $\theta_{ijk}=0$, except either one index $k$ with $\theta_{ijk}=2$, or two distinct indices $k$ with $\theta_{ijk}=1$. Having this in mind, whenever $i\edge j$ is a smoothable edge, we are going to use dark coloring for the angle $i\edge k\edge j$, if $\theta_{ijk}=2$, and light coloring if $\theta_{ijk}=1$; see \autoref{fig:colored_angles_types}.

\begin{example}\label{ex:X4_diagram}
The biresidue matrix on the left corresponds to the smoothing diagram on the right:

$$
\mathcal{X}_4 := 
\begin{pmatrix}
 0 & 2 & -1 & -1 \\
-2 & 0 & 3 & -1 \\
 1 &-3 & 0 & 2 \\
 1 & 1 &-2 & 0
\end{pmatrix}=(b_{ij})_{i,j=0}^3  \hskip0.2\textwidth
\pgfdeclarelayer{background layer}
\pgfdeclarelayer{foreground layer}
\pgfsetlayers{background layer,main,foreground layer}
\begin{tikzpicture}[baseline=-1ex,
rotate=360/2/4,
scale=1.2,line join = round 
    ] 
    
    \begin{pgfonlayer}{main}
\foreach \n in {1,...,4}
{
    \coordinate (v\n) at ({-90+((\n-1)*360/4)}:1.5);
}
\node[below right] at (v1) {$0$};
\node[above right] at (v2) {$1$};
\node[above left] at (v3) {$2$};
\node[below left] at (v4) {$3$};

\foreach \n in {1,...,4}
{
	\foreach \m in {1,...,4}
	\draw (v\n) -- (v\m);
}
\end{pgfonlayer}

\smoothedge{v1}{v2}{v3}{v3} 
\smoothedge{v2}{v3}{v4}{v1}
\smoothedge{v3}{v4}{v2}{v2}
\smoothedge{v4}{v1}{v2}{v3}

\begin{pgfonlayer}{foreground layer}
\foreach \n in {1,...,4}
{
	\draw[fill] (v\n) circle (0.03);
}
\end{pgfonlayer}
 \end{tikzpicture}
$$
\end{example}

The importance of the smoothing diagrams lies in the following result.

\begin{proposition}\label{prop:smoothable_edge_deform}
Let $\X=\mathbb{P}^{\n-1}\times \mathbb{D}^m$, and $\omega_0$ be a semi-toric log symplectic form on $\X$. Let $i\edge j$ be a smoothable edge for $\omega_0$. Then there is a Poisson deformation
$$
\pi(\varepsilon) = \pi_0 + \varepsilon \rho_{ij}, ~~\varepsilon \in \mathbb{C},
$$
where $\pi_0 = \omega_0^{-1}$ and
\begin{equation}\label{eq:smoothable_edge_deform}
\rho_{ij} = \left(\prod_{\ell=1}^m e^{\lambda_\ell x_\ell} \prod_{\substack{k=0\\ k\not=i,j}}^{\n-1} y_k^{\theta_{ijk}} \right) \,\partial_{y_i} \wedge \partial_{y_j},
\end{equation}
for certain constants $\lambda_1,...,\lambda_m\in\mathbb{C}$ (specified in the proof below).
\end{proposition}

\begin{figure}[h]
\centering
\captionsetup{justification=centering}
\begin{subfigure}[t]{0.35\textwidth}
\centering
\begin{tikzpicture}[scale=0.4]
\draw[red, ultra thick] (5,0) -- (-5,0);
\node[right] at (5,0) { $y_i=0$};
\draw[red, ultra thick] (0,5) -- (0,-5);
\node[above] at (0,5) { $y_j=0$};
\end{tikzpicture}
\caption{Polar divisor of $\omega_0$\hspace{0.15\textwidth}~}
\end{subfigure}
\hspace{0.15\textwidth}
\begin{subfigure}[t]{0.35\textwidth}
\centering
\begin{tikzpicture}[scale=0.4]
\draw[ultra thin] (5,0) -- (-5,0);
\node[left] at (-5,0) {~~~~~};
\node[right] at (5,0) { $y_i=0$};
\draw[ultra thin] (0,5) -- (0,-5);
\node[above] at (0,5) { $y_j=0$};
\draw[red, ultra thick, smooth, samples = 200, domain = -5:-1/5] plot(\x, 1/\x);
\draw[red, ultra thick, smooth, samples = 200, domain = 1/5:5] plot(\x, 1/\x);
\end{tikzpicture}
\caption{Polar divisor of $\pi(\varepsilon)^{-1}$, \\ for small $\varepsilon\not=0$}
\end{subfigure}

\caption{Effect of the Poisson deformation in \autoref{prop:smoothable_edge_deform} on the polar divisor \\ near a generic point of $\{y_i=y_j=0\}$.}
\label{fig:smoothing_out_divisor}
\end{figure}
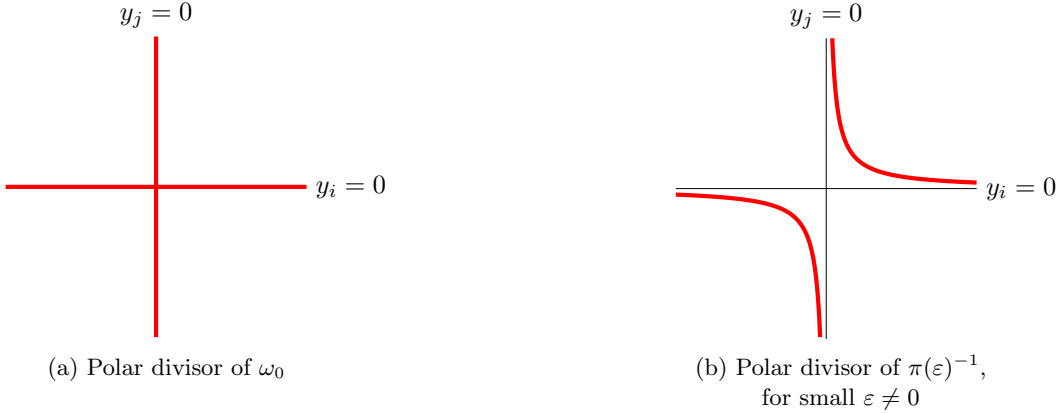

\begin{remark}
The Poisson deformation described in \autoref{prop:smoothable_edge_deform} is smoothing out the irreducible component $\{y_i=y_j=0\}$ of the singular locus of the polar divisor, at least generically (see \autoref{fig:smoothing_out_divisor}). This phenomenon explains the terminology ``smoothable edge'' and ``smoothing diagram''.
\end{remark}

\begin{proof}
The claim can be deduced from \cite{Matviichuk2020}. For reader's convenience, we give a direct argument here.

Without loss of generality, we may assume $i,j\not=0$. Let us use the affine coordinates $z_k = \frac{y_k}{y_0}$, $k=1,...,{\n-1}$. Let $\Pi$ and $\Omega$ be skew-symmetric matrices of size ${\n-1}+m$ such that
$$
\pi_0 = \sum_{1\le k<\ell\le {\n-1}+m} \Pi_{k\ell} ~u_k \wedge u_\ell,~~~
\omega_0 = \sum_{1\le k<\ell\le {\n-1}+m} \Omega_{k\ell}~ u^*_k \wedge u^*_\ell,
$$
where $u_k = z_k\partial_{z_k}$, $u^*_k = d\log(z_k)$, $k=1,...,{\n-1}$, and $u_k = \partial_{x_{{\n-1}+k}}$, $u^*_k = dx_{{\n-1}+k}$, $k=\n,...,{\n-1}+m$. We have $\Omega_{k\ell} = b_{k\ell}$, $1\le k,\ell\le {\n-1}$, and
$$
\rho_{ij} = \left(\prod_{\ell=1}^m e^{\lambda_\ell x_\ell} \prod_{\substack{k=1\\k\not=i,j}}^{\n-1} z_k^{\theta_{ijk}}\right) \partial_{z_i}\wedge \partial_{z_j}.
$$
Let us choose the constants $\lambda_1,..., \lambda_m$ so that
\begin{equation}\label{eq:defn_lambdas}
\begin{pmatrix}
\boldsymbol{\theta} \\
\boldsymbol{\lambda}
\end{pmatrix} = 
\Omega ~\bf{v},
\end{equation}
where $\boldsymbol{\theta}= (\theta_{ij1},...,\theta_{ij,{\n-1}})^T$, $\boldsymbol\lambda = (\lambda_1,...,\lambda_m)^T$, and ${\bf v}=(v_1,...,v_{{\n-1}+m})^T$ with
$$
v_k = \left\{
\begin{matrix}
\frac{1}{b_{ij}}, & \text{ if } k=i, \\
-\frac{1}{b_{ij}}, & \text{ if } k=j,\\
0, & \text{otherwise}.
\end{matrix}
\right.
$$
\autoref{eq:defn_lambdas} implies
$$
\Pi
\begin{pmatrix}
\boldsymbol{\theta} \\
\boldsymbol{\lambda}
\end{pmatrix} = 
\bf{v},
$$
which can be rewritten in terms of the Schouten bracket as follows
$$
\big[\pi_0, \rho_{ij}\big] = \left(\sum_{k=1}^{{\n-1}+m} v_k \,u_k\right)\wedge \rho_{ij}=\dfrac{1}{b_{ij}} \left(z_i\partial_{z_i} - z_j\partial_{z_j}  \right) \wedge \rho_{ij} = 0.
$$
It remains to note that $[\rho_{ij},\rho_{ij}]=0$.
\end{proof}

\begin{example}\label{ex:X4_first_order}
Consider the following log symplectic form on $\mathbb{P}^3 \times \mathbb{D}^1$
$$
\omega_0 = \dfrac{1}{12} \left(\sum_{1\le i<j\le 3} b_{ij} \dfrac{dz_i}{z_i} \wedge \dfrac{dz_j}{z_j}\right) + \dfrac{1}{24} \left( -3 \dfrac{dz_1}{z_1} + 3 \dfrac{dz_2}{z_2} - \dfrac{dz_3}{z_3}\right) \wedge dx,
$$
where the biresidues $b_{ij}$ are as in \autoref{ex:X4_diagram}, $z_i=\frac{y_i}{y_0}$, $i=1,2,3$, are affine coordinates on $\mathbb{P}^3$, and $x$ is the standard coordinate on $\mathbb{D}^1$. Let us invert $\omega_0$ to obtain
$$
\pi_0 = - 2 \left(z_1\partial_{z_1} \wedge z_2 \partial_{z_2} + 3~z_1\partial_{z_1}\wedge z_3 \partial_{z_3} + 3~z_2 \partial_{z_2} \wedge z_3 \partial_{z_3}\right) + 4 \left(2~z_1 \partial_{z_1} + z_2 \partial_{z_2} + 3~z_3\partial_{z_3} \right) \wedge \partial_x.
$$
Using the homogeneous coordinates $y_i$, $i=0,1,2,3$, on $\mathbb{P}^3$ we get
$$
\pi_0 = - \Big(2~y_0 \partial_{y_0} \wedge y_1 \partial_{y_1} + y_0 \partial_{y_0} \wedge y_2 \partial_{y_2}  - 3~y_0 \partial_{y_0} \wedge y_3 \partial_{y_3} + y_1 \partial_{y_1} \wedge y_2 \partial_{y_2} + y_1 \partial_{y_1} \wedge y_3 \partial_{y_3} + 2~y_2 \partial_{y_2} \wedge y_3 \partial_{y_3}\Big) -
$$
$$
-2~ \Big( 3 ~y_0 \partial_{y_0} - y_1\partial_{y_1} + y_2 \partial_{y_2} -3~y_3\partial_{y_3}\Big) \wedge \partial_x.
$$

As we have seen in \autoref{ex:X4_diagram}, the log symplectic form $\omega_0$ has four smoothable edges $i\edge (i+1)\text{~mod~} 4$, $i=0,1,2,3$. The corresponding Poisson deformations are given by
$$
\begin{matrix}
     \rho_{01} &=& e^{-x} ~y_2^2~ \partial_{y_0} \wedge \partial_{y_1},  \\
     \rho_{12} &=& e^{x}~ y_3 y_0 ~\partial_{y_1} \wedge \partial_{y_2}, \\
     \rho_{23} &=& e^{-x}~  y_1^2 ~\partial_{y_2} \wedge \partial_{y_3}, \\
     \rho_{30} &=& e^{x}~ y_1 y_2 ~\partial_{y_3} \wedge \partial_{y_0}.
\end{matrix}
$$
\end{example}

Poisson deformations arising from smoothable edges exhibit a striking property -- they are jointly unobstructed. This claim, unlike \autoref{prop:smoothable_edge_deform}, does not seem to admit an elementary proof. The proof we know uses in an essential way the structure of $L_\infty$ algebra that the Schouten bracket induces on the Poisson cohomology via homotopy transfer.

\begin{theorem}[\cite{Matviichuk2020}]\label{thm:deforming_several_smootables}
Let $\X=\mathbb{P}^{\n-1}\times \mathbb{D}^m$, and $\omega_0$ be a semi-toric log symplectic form on $\X$.
Let the edges $i_1\edge j_1$, ..., $i_s\edge j_s$ be smoothable for $\omega_0$. Let $\pi_0=\omega_0^{-1}$ and 
$
\pi_1 = \sum_{a=1}^s \mu_{i_aj_a} ~\rho_{i_aj_a}
$, where each $\rho_{i_aj_a}$ is the bivector defined in \eqref{eq:smoothable_edge_deform}, and $\mu_{i_aj_a}\in\mathbb{C}^*$.

Then there exists a formal Poisson deformation 
\begin{equation}\label{eq:deforming_several_smoothables}
\pi(\varepsilon)=\sum_{k=0}^\infty \varepsilon^k \pi_k,
\end{equation} 
where each $\pi_k$ is a holomorphic bivector on $\X$.
\end{theorem}

\begin{proof}
The claim follows from \cite[Corollary 6.23]{Matviichuk2020}. The condition 2 of the corollary is satisfied by the argument given in Example 6.24. Corollary 6.23 implies that the deformation functor governing the Poisson deformations of $(\X, \pi_0)$ admits the constructible vector bundle description given in Theorem 6.16. The Poisson deformation in \eqref{eq:deforming_several_smoothables} corresponds to the central fiber of this bundle.
\end{proof}

\begin{conjecture}
The power series in \eqref{eq:deforming_several_smoothables} can be chosen to be convergent, at least after shrinking the polydisc.
\end{conjecture}

For the case $m=0$, one can prove the conjecture using compactness of $\X$ and the Artin approximation theorem \cite{Artin68}. Additionally, we verify the conjecture for several examples with $m\ge1$ below.

\begin{remark}
The methods described \cite{Matviichuk2020} do not immediately provide an explicit formula for the Poisson deformation \eqref{eq:deforming_several_smoothables}. We hope that a more careful analysis of the $L_\infty$ algebra on the Poisson cohomology will eventually lead to such a formula. Until then, one must determine the power series on a case-by-case basis.
\end{remark}


\begin{example} \label{ex:X4_full_deform}
Continuing with \autoref{ex:X4_first_order}, 
let $\pi_1= \sum_{i=0}^3 \rho_{i,i+1}$. Let us look  for a Poisson deformation $\pi(\varepsilon) = \pi_0 + \varepsilon \pi_1 + \varepsilon^2 \pi_2 +... $ satisfying the following ansatz
$$
\pi(\varepsilon) = \pi_0 + \sum_{i=0}^3 \xi_i(x,\varepsilon) \rho_{i,i+1} + \eta(x,\varepsilon) \sum_{i=0}^3 y_i\partial_{y_i} \wedge y_{i+1} \partial_{y_{i+1}}.
$$
\end{example}
After working out the Schouten integrability $[\pi(\varepsilon),\pi(\varepsilon)]=0$, one finds that the deformation $\pi(\varepsilon)$ is Poisson if and only if $\xi_0,\xi_1,\xi_2,\xi_3,\eta$ satisfy the following system of ODEs (all derivatives are with respect to $x$):
\begin{equation*}
\begin{aligned}
    \eta' &= - \frac{1}{8}e^{2x}\xi_1 \xi_3, \\
    \xi_0' &= \frac{1}{2}  \xi_0 \eta, \\
    \xi_1'& = -\frac{1}{3}  \xi_1\eta, \\
    \xi_2' &= \frac{1}{2}  \xi_2\eta, \\
    \xi_3'& = -  \xi_3 \eta+ e^{-3x}\xi_0 \xi_2. \\
\end{aligned} 
\end{equation*}

The Picard-Lindel\"{o}f theorem implies that one can solve these ODEs uniquely, after specifying the initial conditions at $x=0$, and possibly shrinking the disc $\mathbb{D}^1$. For instance, we can specify $\xi_i\big|_{x=0} = \varepsilon$, $\eta\big|_{x=0} = 0$.

To find where the Poisson structure $\pi=\pi(\varepsilon)$ on $\X=\mathbb{P}^3\times \mathbb{D}^1$ drops rank, we calculate its Pfaffian

\begin{wrapfigure}{r}{0.35\textwidth}
\vskip-0.2cm
\includegraphics[scale=0.7]{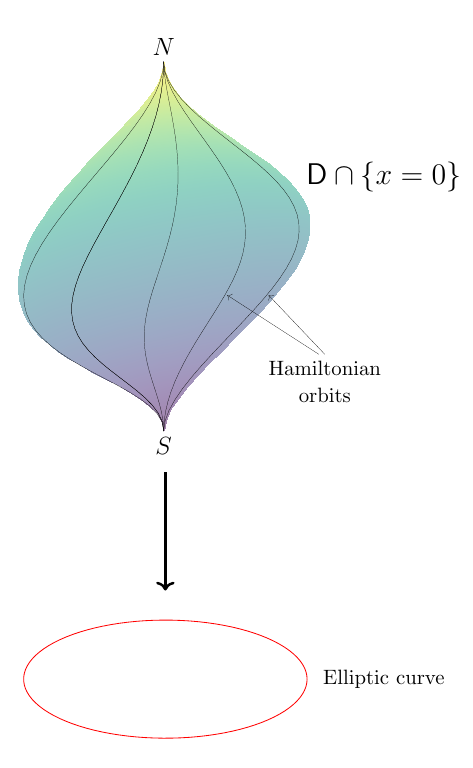}
%
%
%
%
%
%
%
%
%
\caption{Degeneracy divisor of the Poisson deformation $\pi(\varepsilon)$}
\label{fig:X4Divisor}
\end{wrapfigure}

$$
\displaystyle \pi \wedge \pi\wedge \sum_{i=0}^3 y_i\partial_{y_i} = {\rm Pf}(\pi) \bigwedge_{i=0}^3 \partial_{y_i} \wedge \partial_x,
$$
$$
{\rm Pf}(\pi) 
= 32e^x(-3+2\eta) y_0y_1y_2y_3 + 16 \xi_0 y_2^3y_3 + 
$$
$$
+24 e^{2x} \xi_1 y_0^2 y_3^2 + 16 \xi_2 y_0y_1^3 +8 e^{2x} \xi_3 y_1^2 y_2^2.
$$

Applying the Jacobian criterion to $\D = \{{\rm Pf}(\pi)=0\}$, we obtain that for small values of $\varepsilon\not=0$ and $x$, the divisor $\D$ is smooth away from $S\times \mathbb{D}^1$ and $N \times \mathbb{D}^1$, for $S=[1:0:0:0]$, $N=[0:0:0:1]$, and likewise, the quartic surface $\mathsf{Q} = \D\cap\{x=0\}\subset\mathbb{P}^3$ is smooth away from $S,N$. One can check that $\mathsf{Q}$ has an elliptic $\widetilde{E}_8$  singularity (studied in \cite{Etingof2010,Pym2017}) at each of these points. The space of codimension $2$ symplectic leaves is isomorphic to the quotient of the smooth surface $\mathsf{Q}^\circ=\mathsf{Q} \setminus \{S,N\}$ by the orbits of the Hamiltonian vector field 
$$
u = [\pi, x] = [\pi_0,x] = -2 \Big( 3 ~y_0 \partial_{y_0} - y_1\partial_{y_1} + y_2 \partial_{y_2} -3~y_3\partial_{y_3}\Big).
$$
The quotient is Hausdorff, since the orbits are closed in $\mathsf{Q}^\circ$. It is also compact, has dimension one and admits a non-vanishing vector field, namely the modular vector field of $\pi$. Therefore, the space of codimension $2$ leaves of $\pi$ is an {\bf elliptic curve}. We illustrate this on \autoref{fig:X4Divisor}, where the surface $\mathsf{Q}=\D\cap\{x=0\}\subset\mathbb{P}^3$ takes on the shape of an onion.


\newpage
\section{Smoothable cycles}\label{sec:smoothableCycles}

The following section is entirely elementary, and deals with the combinatorics of smoothing diagrams. By a biresidue matrix we mean a skew-symmetric  matrix with complex entries, referred to as biresidues.

\begin{lemma}[\cite{Matviichuk2020}]\label{lm:pos_rat_valency}
Let $B=(b_{ij})_{i,j=0}^{\n-1}$ be a biresidue matrix. Then 
\begin{enumerate}
\item For each vertex $j$, there are at most two smoothable edges containing $j$.
\item If $i\edge j$ and $j\edge k$ are two smoothable edges for $B$, then $b_{ij}$ is a positive rational multiple of $b_{jk}$.
\end{enumerate}
\end{lemma}

\begin{proof}
Let us prove 2 first. Suppose that
$$
\dfrac{b_{jk}+b_{ki}}{b_{ij}}=p\in\mathbb{Z}_{\ge0} ~~~~\text{and} ~~~~\dfrac{b_{ki}+b_{ij}}{b_{jk}}=q\in\mathbb{Z}_{\ge0},
$$
By subtracting two equalities $b_{jk}+b_{ki} = p~b_{ij}$ and $b_{ki}+b_{ij} = q~b_{jk}$, we obtain $(1+q) ~b_{jk} = (1+p) ~b_{ij}$.

To prove 1, let us assume that $j\edge i$, $j\edge \ell$, $j\edge k$ are all different smoothable edges. Then, by 2, we get $b_{ji}$ and $b_{j\ell}$ are negative multiples of each other, and so are $b_{j\ell}$ and $b_{jk}$. However, those two facts together imply that $b_{ji}$ and $b_{jk}$ are positive multiples of each other, which contradicts 2.
\end{proof}

\autoref{lm:pos_rat_valency} part 1 implies that each connected component of the graph formed by the smoothable edges is either a chain, or a cycle. We are going to call a cycle of smoothable edges a \textit{smoothable cycle}. If the smoothable cycle involves all the vertices $0,1,...,\n-1$, we are going to say that the biresidue matrix, or the corresponding smoothing diagram, has a \textit{full smoothable cycle}. In what follows we are going to classify the biresidue matrices with zero sum in each row that have a full smoothable cycle (as announced in \autoref{thm:smoothableCycles}). Let us start by listing all such matrices, and then proceed to proving that the list is exhaustive.

Everywhere below we will number the rows and columns of a matrix starting from zero, rather than one. 
For a row vector ${\bf r}=(r_0,r_1,...,r_{\n-1})$, we denote by ${\rm T}\,{\bf r}$, or ${\rm T}({\bf r})$, the cyclic shift $(r_{\n-1},r_0,r_1,...,r_{\n-2})$ of ${\bf r}$ by one position to the right.

\noindent {\bf Family $\mathcal{C}_{\n,k}$:} For any numbers $\n\ge 3$ and $1\le k < \n/2$, we define the matrix $\mathcal{C}_{\n,k}$ whose zeroth row is
\begin{equation}\label{eq:row_basic}
{\bf r}(k) := (0 ~\underbrace{1~ ... ~1}_{k} ~ \underbrace{0~ ... ~0}_{\n-1-2k} ~\underbrace{-1~...~-1}_{k} ),
\end{equation}
and whose $i$-th row is ${\rm T}^i {\bf r}(k)$. For instance,

$$
\mathcal{C}_{5,1} = 
\begin{pmatrix}
 0 & 1 & 0 & 0 & -1 \\
 -1 & 0 & 1 & 0 & 0 \\
 0 & -1 & 0 & 1 & 0\\
 0 & 0 & -1 & 0 & 1 \\
 1 & 0 & 0 & -1 & 0 \\
\end{pmatrix},
\hskip2cm
\mathcal{C}_{7,3} = 
\begin{pmatrix}
 0 & 1 & 1 & 1 & -1 & -1 & -1 \\
 -1 & 0 & 1 & 1 & 1 & -1 & -1 \\
 -1 & -1 & 0 & 1 & 1 & 1 & -1\\
 -1 & -1 & -1 & 0 & 1 & 1 & 1\\
 1 & -1 & -1 & -1 & 0 & 1 & 1\\
 1 & 1 & -1 & -1 & -1 & 0 & 1\\
 1 &1 & 1 & -1 & -1 & -1 & 0\\
\end{pmatrix}.
$$

\noindent{\bf Family $\mathcal{C}_{\n,k,I}$:} For $\n\ge 3$, $1\le k <\n/2$ and $I=\{i_0,i_1,...,i_{d-1}\}$ being a list of $0$'s and $1$'s of length $d = \gcd(\n,k)$, we define the matrix $\mathcal{C}_{\n,k,I}$ as follows. The first $d$ rows of this matrix are

	$$
	{\bf r}_0={\bf r}(k-i_0), {\bf r}_1 ={\rm T} ({\bf r}(k-i_1)),  ... , {\bf r}_{d-1} = {\rm T}^{d-1}({\bf r}(k-i_{d-1}))
	$$
	where ${\bf r}(l)$ is defined as in \eqref{eq:row_basic}. For $d\le i\le \n-1$, the $i$-th row ${\bf r}_i$ of $\mathcal{C}_{\n, k,I}$ is defined using the following $d$-periodicity condition: ${\bf r}_i = {\rm T}^d {\bf r}_{i-d}$.
	
\noindent For example, $\mathcal{C}_{\n,k,\{0,0,...,0\}} =\mathcal{C}_{\n,k}$; $\mathcal{C}_{\n,k,\{1,1,...,1\}}=\mathcal{C}_{\n,k-1}$;
	$$
	\mathcal{C}_{9,3, \{0,1,1\}} = 
	\begin{pmatrix}
	0  & 1 & 1 & 1 & 0 & 0 &-1 &-1 &-1 \\
	-1 & 0 & 1 & 1 & 0 & 0 & 0 & 0 &-1 \\
	-1 &-1 & 0 & 1 & 1 & 0 & 0 & 0 & 0 \\ 
	-1 &-1 &-1 & 0 & 1 & 1 & 1 & 0 & 0 \\ 
	 0 & 0 &-1 &-1 & 0 & 1 & 1 & 0 & 0 \\ 
	 0 & 0 & 0 &-1 &-1 & 0 & 1 & 1 & 0 \\ 
	 1 & 0 & 0 &-1 &-1 &-1 & 0 & 1 & 1 \\ 
	 1 & 0 & 0 & 0 & 0 &-1 &-1 & 0 & 1 \\ 
	 1 & 1 & 0 & 0 & 0 & 0 &-1 &-1 & 0 \\ 
	\end{pmatrix}
	=
	\begin{pmatrix}
	{\bf r}(3) \\
	{\rm T}({\bf r}(2)) \\
	{\rm T}^2({\bf r}(2)\\ 
	{\rm T}^3({\bf r}(3)) \\ 
	{\rm T}^4({\bf r}(2))\\ 
	{\rm T}^5({\bf r}(2))\\ 
	{\rm T}^6({\bf r}(3))\\ 
	{\rm T}^7({\bf r}(2))\\ 
	{\rm T}^8({\bf r}(2)) \\ 
	\end{pmatrix}.
	$$

\noindent{\bf Sporadic example $\mathcal{X}_4$} was presented in \autoref{ex:X4_diagram}.

\noindent{\bf Sporadic example $\mathcal{X}_5$} is defined as follows
$$
\mathcal{X}_5 = 
\begin{pmatrix}
0 & 1 & 1 & -1 & -1 \\
-1& 0 & 2 & 0 &-1 \\
-1&-2 & 0 & 2 & 1 \\
 1& 0 &-2 & 0 & 1 \\
 1 & 1 &-1&-1&0
\end{pmatrix}.
$$

\noindent{\bf Family $\mathcal{Y}_{2(2k+1)}$:}~
For $k\ge1$, the first two rows ${\bf r}_0$, ${\bf r}_1$ of $\mathcal{Y}_{2(2k+1)}$ are defined follows:

$$
\begin{pmatrix}
 0 & 2 & 1 & 1 & .... & 1 &-1 &-1 &-1& ... & -1 \\
-2 & 0 & 1 & 1 & .... & 1 & 1 & 1 &-1& ... & -1
\end{pmatrix},
$$
where the zeroth row has $\n/2-2$ ones, for $\n=2(2k+1)$. For $2\le i\le \n-1$, the $i$-th row of the matrix $\mathcal{Y}_{2(2k+1)}$ is determined by the following periodicity condition: ${\bf r}_i = {\rm T}^2({\bf r}_{i-2})$. For instance:
$$
\mathcal{Y}_6 = 
\begin{pmatrix}
 0 & 2 & 1 &-1 &-1 &-1 \\
-2 & 0 & 1 & 1 & 1 &-1 \\
-1 &-1 & 0 & 2 & 1 &-1 \\
 1 &-1 &-2 & 0 & 1 & 1 \\
 1 &-1 &-1 &-1 & 0 & 2 \\
 1 & 1 & 1 &-1 &-2 & 0 \\
\end{pmatrix},
~~~~~~~~
\mathcal{Y}_{10} = 
\begin{pmatrix}
 0 & 2 & 1 & 1 & 1 &-1 &-1 &-1 &-1 &-1 \\
-2 & 0 & 1 & 1 & 1 & 1 & 1 &-1 &-1 &-1 \\
-1 &-1 & 0 & 2 & 1 & 1 & 1 &-1 &-1 &-1 \\
-1 &-1 &-2 & 0 & 1 & 1 & 1 & 1 & 1 &-1 \\
-1 &-1 &-1 &-1 & 0 & 2 & 1 & 1 & 1 &-1 \\
 1 &-1 &-1 &-1 &-2 & 0 & 1 & 1 & 1 & 1 \\
 1 &-1 &-1 &-1 &-1 &-1 & 0 & 2 & 1 & 1 \\
 1 & 1 & 1 &-1 &-1 &-1 &-2 & 0 & 1 & 1 \\
 1 & 1 & 1 &-1 &-1 &-1 &-1 &-1 & 0 & 2 \\
 1 & 1 & 1 & 1 & 1 &-1 &-1 &-1 &-2 & 0 \\
\end{pmatrix}.
$$

\noindent{\bf Family $\mathcal{Z}_{5(2k+1)}$:} For $k\ge 1$, the first five rows of $\mathcal{Z}_{5(2k+1)}$ are as follows:
$$
\begin{pmatrix}
{\bf r}_0 \\
{\bf r}_1 \\
{\bf r}_2 \\
{\bf r}_3 \\
{\bf r}_4 \\
\end{pmatrix}=
\begin{pmatrix}
0 & 1 & 1 & 1 & 1 & 1&...& 1 & -1 & -1& -1& -1&-1&-1&...& -1  \\
-1& 0 & 2 & 2 & 1 & 1&...& 1 & -1 & -1& -1& -1&-1&-1&...& -1 \\
-1&-2 & 0 & 2 & 1 & 1&...& 1 &  1 &  1& -1& -1&-1&-1&...&-1\\
-1&-2 &-2 & 0 & 1 & 1&...& 1 &  1 &  1&  1&  1&-1&-1&...&-1\\
-1 & -1 &-1&-1& 0 & 1&...& 1 &  1 &  1&  1& 1&-1&-1&...&-1\\
\end{pmatrix},
$$
where the zeroth row has $(\n-1)/2$ ones, for $\n=5(2k+1)$. For $5\le i\le \n-1$, the $i$-th row of the matrix $\mathcal{Z}_{5(2k+1)}$ is determined by the following periodicity condition: ${\bf r}_{i}={\rm T}^5({\bf r}_{i-5})$. For instance:
$$
\mathcal{Z}_{15} = 
\begin{pmatrix}
0  & 1 & 1 & 1 & 1 & 1 & 1 & 1 & -1 & -1 &-1 &-1 &-1 &-1 &-1 \\
-1 & 0  & 2 & 2 & 1 & 1 & 1 & 1 & -1 & -1 & -1 &-1 &-1 &-1 &-1 \\
-1 &-2 & 0  & 2 & 1 & 1 & 1 & 1 & 1 & 1 & -1 & -1 &-1 &-1 &-1  \\
-1 & -2 &-2 & 0  & 1 & 1 & 1 & 1 & 1 & 1 & 1 & 1 & -1 &-1 &-1  \\
-1 & -1 & -1 &-1 & 0  & 1 & 1 & 1 & 1 & 1 & 1 & 1 & -1 & -1 &-1  \\
-1 &-1 & -1 & -1 &-1 & 0  & 1 & 1 & 1 & 1 & 1 & 1 & 1 & -1 & -1  \\
-1 & -1 &-1 & -1 & -1 &-1 & 0  & 2 & 2 & 1 & 1 & 1 & 1 & -1 & -1   \\
-1 & -1 & -1 &-1 & -1 & -1 &-2 & 0  & 2 & 1 & 1 & 1 & 1 & 1 & 1   \\
1 & 1 & -1 & -1 &-1 & -1 & -2 &-2 & 0  & 1 & 1 & 1 & 1 & 1 & 1   \\
1 & 1 & -1 & -1 & -1 &-1 & -1 & -1 &-1 & 0  & 1 & 1 & 1 & 1 & 1    \\
1 & 1 & 1 & -1 & -1 & -1 &-1 & -1 & -1 &-1 & 0  & 1 & 1 & 1 & 1    \\
1 & 1 & 1 & -1 & -1 & -1 & -1 &-1 & -1 & -1 &-1 & 0  & 2 & 2 & 1    \\
1 & 1 & 1 & 1 & 1 & -1 & -1 &-1 &-1 & -1 & -1 &-2 & 0  & 2 & 1    \\
1 & 1 & 1 & 1 & 1 & 1 & 1 & -1 &-1 &-1 & -1 & -2 &-2 & 0  & 1   \\
1 & 1 & 1 & 1 & 1 & 1 & 1 & -1 & -1 &-1 &-1 & -1 & -1 &-1 & 0  \\
\end{pmatrix}.
$$

We will say that two biresidue matrices $B$ and $B'$ are \textit{equivalent}, if there is a permutation matrix $P$ such that $B'=P^{-1}BP$; and \textit{projectively equivalent} if there is a permutation matrix $P$ and a scalar $\lambda\in \mathbb{C}\setminus\{0\}$ such that $B'=\lambda P^{-1}BP$.

\newtheorem*{theoremC}{\autoref{thm:smoothableCycles}}
\begin{theoremC}
\textit{ Let $B$ be a biresidue matrix whose rows sum to zero. Then $B$ has a full smoothable cycle if and only if it is projectively equivalent to one of the following matrices:
\item[-] $\mathcal{C}_{\n,k,I}$, where $1\le k< \n/2$, and $I$ is a sequence of zeros and ones of length $d=\gcd(\n,k)$,
\item[-] $\mathcal{X}_4$, $\mathcal{X}_5$, 
\item[-] $\mathcal{Y}_{2(2k+1)}$, $k\ge 1$, 
\item[-] $\mathcal{Z}_{5(2k+1)}$, $k\ge 1$.}
\end{theoremC}

%

\begin{proof}[Proof of \autoref{thm:smoothableCycles}]
For sufficiency, it is fairly straightforward to verify that each biresidue matrix in the list possesses a cycle comprising the smoothable edges $i\edge (i+1)\mod n$, for $i=0,1,...,n-1$. Refer to \autoref{fig:smoothing_diagrams_gallery} for examples of the associated smoothing diagrams. Now, let's move on to proving necessity.

Let $B=(b_{ij})_{i,j=0}^{\n-1}$ be a biresidue matrix with zero row sum admitting a full smoothable cycle. After conjugating it via a permutation matrix, we may assume that the smoothable cycle is formed by the edges $i\edge (i+1)\mod{\n}$, $i=0,1,...,\n-1$. \autoref{lm:pos_rat_valency}, part 2, implies that the entries $b_{i,i+1}$, $i=0,1,...,\n-1$, are all rational multiples of each other, where here and everywhere below we consider the indices modulo $\n$. By rescaling, we may assume that $b_{i,i+1}>0$ for each $0\le i\le \n-1$. For each $i,j$, smoothability of the edge $j\edge (j+1)$ implies that $b_{i,j}-b_{i,j+1} \in \{0,b_{j,j+1},2b_{j,j+1}\} \subset \mathbb{R}_{\ge0}$. Therefore, for each $i$, the sequence $b_{i,\bullet}=\{b_{i,i+1},b_{i,i+2},...,b_{i,i-1}\}$ is non-increasing. Also, we have $\sum_{j=0}^{\n-1} b_{i,j}=0$ for each $i$. To carry out the proof, we will consider different cases. The cases and their outcomes can be summarized as follows

\begin{table}[h]
    \centering
    \begin{tabular}{|c|c|c|c|c|}
    \hline
        Case 1& \multicolumn{2}{|c|}{Case 2} &\multicolumn{2}{|c|}{Case 3} \\
        $\forall i: b_{i,i+1}=b_{i+1,i+2}$& \multicolumn{2}{|c|}{each vertex has two colored angles} & \multicolumn{2}{|c|}{there is a vertex $i$ with one colored angle} \\
        \hline
        \multirow{4}{2em}{$\mathcal{C}_{n,k,I}$} & Subcase 2.1 &  Subcase 2.2 & Subcase 3.1 & Subcase 3.2 \\
            & $\exists i :b_{i,i+1}=b_{i+1,i+2}$ & $\forall i : b_{i,i+1}\not=b_{i+1,i+2}$&$b_{i,i+1}\not= -b_{i,i-1}$ & $b_{i,i+1}=-b_{i,i-1}$ \\ 
            \cline{2-5} &&&&\\[-1em]
      &  $\mathcal{C}_{n,k,I}$ & $\mathcal{Y}_{2(2k+1)}$ & $\mathcal{X}_4$ &$\mathcal{C}_{2k+1,k},\mathcal{X}_5, \mathcal{Z}_{5(2k+1)}$\\
      \hline
        
    \end{tabular}
    \label{tab:my_label}
\end{table}

\begin{figure}
\centering
\begin{subfigure}{0.3\textwidth}
\sageNgon[1.3]{6}{
\smoothedge{v1}{v2}{v4}{v5} 
\smoothedge{v2}{v3}{v4}{v1}
\smoothedge{v3}{v4}{v6}{v1}
\smoothedge{v4}{v5}{v6}{v3}
\smoothedge{v5}{v6}{v2}{v3}
\smoothedge{v6}{v1}{v2}{v5}
}
\caption{$\mathcal{Y}_{6}$~~~~~~~~~~}
\end{subfigure}
\begin{subfigure}{0.3\textwidth}
\sageNgon[1.5]{7}{
\smoothedge{v1}{v2}{v5}{v5} 
\smoothedge{v2}{v3}{v6}{v6}
\smoothedge{v3}{v4}{v7}{v7}
\smoothedge{v4}{v5}{v1}{v1}
\smoothedge{v5}{v6}{v2}{v2}
\smoothedge{v6}{v7}{v3}{v3}
\smoothedge{v7}{v1}{v4}{v4}
}
\caption{$\mathcal{C}_{7,3}$}
\end{subfigure}
\begin{subfigure}{0.3\textwidth}
\sageNgon[1.9]{9}{
\smoothedge{v1}{v2}{v7}{v8} 
\smoothedge{v2}{v3}{v5}{v9}
\smoothedge{v3}{v4}{v6}{v7}
\smoothedge{v4}{v5}{v1}{v2}
\smoothedge{v5}{v6}{v8}{v3}
\smoothedge{v6}{v7}{v9}{v1}
\smoothedge{v7}{v8}{v4}{v5}
\smoothedge{v8}{v9}{v2}{v6}
\smoothedge{v9}{v1}{v3}{v4}
}
\caption{$\mathcal{C}_{9,3,\{0,1,1\}}$}
\end{subfigure}\\
\vskip1cm
\begin{subfigure}{0.75\textwidth}
\sageNgon[4]{15}{
\smoothedge{v1}{v2}{v3}{v4}
\smoothedge{v2}{v3}{v9}{v10}
\smoothedge{v3}{v4}{v11}{v12}
\smoothedge{v4}{v5}{v2}{v3}
\smoothedge{v5}{v6}{v13}{v13}
\smoothedge{v6}{v7}{v8}{v9}
\smoothedge{v7}{v8}{v14}{v15}
\smoothedge{v8}{v9}{v1}{v2}
\smoothedge{v9}{v10}{v7}{v8}
\smoothedge{v10}{v11}{v3}{v3}
\smoothedge{v11}{v12}{v13}{v14}
\smoothedge{v12}{v13}{v4}{v5}
\smoothedge{v13}{v14}{v6}{v7}
\smoothedge{v14}{v15}{v12}{v13}
\smoothedge{v15}{v1}{v8}{v8}
}
\caption{$\mathcal{Z}_{15}$}
\end{subfigure}
\centering

\caption{ {Examples of smoothing diagrams having full smoothable cycles. \\ The vertices are labeled counterclockwise, the bottom horizontal edge being $(\n-1) \,\includegraphics{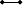}\, 0$.}}
\label{fig:smoothing_diagrams_gallery}
\end{figure}

Cases 2 and 3 together cover all the possibilities, while Case 1 is auxiliary. It is considered separately, since some of the Subcases in Cases 2, 3 reduce to Case 1.

\textbf{Case 1.} The equality $b_{i,i+1}=b_{i+1,i+2}$ holds for each $i$.

By rescaling all the entries of $B$ simultaneously, we can assume that $b_{i,i+1}=1$ for each $i$.

For each $i,j$, smoothability of the edge $j \edge (j+1)$ implies that $b_{i,j}-b_{i,j+1} \in \{0,b_{j,j+1},2b_{j,j+1}\} =\{0,1,2\}$. For each $i$, the biresidues $b_{i,i+1}, b_{i,i+2}, ..., b_{i,i-1}$ form a non-increasing sequence that starts with $1$, ends with $-1$ and each $b_{i,j}$ is either $1$, or $-1$, or $0$. For each $i$, let us denote by $h(i)$ the number of ones in this sequence, i.e. the largest $1\le h < \n$ such that $b_{i,i+h}=1$. Note that since $\sum_{j=0}^{\n-1} b_{i,j}=0$, the number of $-1$'s in the sequence $b_{i,i+1}, b_{i,i+2}, ..., b_{i,i-1}$ also equals $h(i)$.

We claim that for each $i$, one has

1) $h(i+1)-h(i)\in \{-1,0,1\}$.

2) $h(i+ h(i))\ge h(i)$ and $h(i- h(i))\ge h(i)$.

Proof of 1): smoothability of the edge $i\edge (i+1)$ implies that the number of $1$'s in the $(i+1)$-th row is at least $h(i)-1$. The fact that the smoothable edge $i\edge (i+1)$ has at most two colored angles implies that $h(i+1)\le h(i)+1$.

Proof of 2): since $b_{i,i+h(i)}=1$, by skew symmetry we have $b_{i+h(i),i}=-1$. By monotony of $b_{i,\bullet}$ we get that $b_{i+h(i),j}=-1$ for $j=i,i+1,i+2, ..., i+h(i)-1$. Therefore, $h(i+h(i))\ge h(i)$. The inequality $h(i-h(i))\ge h(i)$ is proved analogously.

Let us choose an index $i_0$ such that $h(i_0)=\max_{0\le i <\n} h(i)$. Denoting $k=h(i_0)$, we get from the claim 2) above that the set of indices $i$ with $h(i)=k$ is $k$-periodic. Combining this with $\n$-periodicity of the indices, we get that this set is $d$-periodic, for $d=\gcd(\n,k)$. In particular, if $k$ is relatively prime with $\n$, then we obtain that the biresidue matrix $B$ equals $\mathcal{C}_{\n,k}=\mathcal{C}_{\n,k,\{0\}}$. 

Let us now consider the situation when $d>1$. Claim 1) above implies that for each $i\equiv (i_0\pm 1) {\rm~mod~} d$ one has $h(i)\ge k-1$. We further claim that for each $i\equiv (i_0\pm 2) {\rm~mod~} d$ one has $h(i)\ge k-1$. Indeed, if $h(i_0+1)=k$, then $h(i_0+2)\ge k-1$ by 1). If $h(i_0+1)=k-1$, then $h(i_0+1+k)=k-1$, and so by 2) we have $h(i_0+2)= h(i_0+1+k-h(i_0+1+k))\ge h(i_0+1+k)=k-1$. The inequality $h(i_0-2)\ge k-1$ is proved analogously. The inequalities $h(i_0\pm 2)\ge k-1$ imply then $h(i)\ge k-1$ for any $i\equiv (i_0\pm 2){\rm~mod~}d$ by 2). Inductively, we prove that $h(i)\ge k-1$ for all $i$. Since the set of indices $i$ with $h(i)=k$ is $d$-periodic, the set of $i$ with $h(i)=k-1$ is $d$-periodic, too. Therefore, the biresidue matrix $B$ is equal to $\mathcal{C}_{\n,k,\{s_0,...,s_{d-1}\}}$, where $s_l=0$ if $h(i_0+l)=k$ and $s_l=1$ if $h(i_0+l)=k-1$.

This finishes the proof in Case 1.

\textbf{Case 2.} Each vertex $v$ has exactly two distinct colored angles of the form $i\edge v\edge j$. 

Note that, since there are overall $2\n$ colored angles (counting with multiplicity) and $\n$ vertices, the assumption of Case 2 implies that there are only lightly colored angles. Another consequence of the assumption of Case 2 is that for each $i$ the sequence $b_{i,\bullet}=\{b_{i,i+1}, b_{i,i+2}, ..., b_{i,i-1}\}$ assumes exactly three values.

\textbf{Subcase 2.1.} Additionally assume that there is an $i$ such that $b_{i,i+1}=b_{i+1,i+2}$. 

Without loss of generality we assume that $b_{i,i+1}=1$. Then the sequence $b_{i+1,i+2}, b_{i+1,i+3}, ... , b_{i+1,i}$ starts with $1$, ends with $-1$, is non-increasing and one has $b_{i+1,j}\not=b_{i+1,j+1}$ for exactly two values of $j$ (for those $j$ that $j\edge (i+1)\edge (j+1)$ is a colored angle). Therefore, the elements of this sequence assume exactly three values, two of which are $1$ and $-1$. 

We claim that the third value occurring in $b_{i+1,\bullet}$ has to be $0$. On the contrary, suppose the third value is $w\in (-1,1)\setminus \{0\}$. The sequence $b_{i,\bullet} =\{b_{i,i+1}, b_{i,i+2}, ... , b_{i,i-1}\}$ assumes exactly three values, one of which is $b_{i,i+1}=1$. Moreover, smoothability of the edge $i-(i+1)$ implies $b_{i+1,j}-b_{i,j}\in \{1,2\}$ for two indices $j_1,j_2\in \{i+2,i+3,...,i-1\}$ and $b_{i+1,j}-b_{i,j}=0$ for remaining $j\in \{i+2,i+3,...,i-1\}$. First, note that $j_1$ cannot be equal to $j_2$, because otherwise we would have $b_{i+1,j_1}-b_{i,j_1}=2\implies b_{i,j_1}= -1, b_{i+1,j_1}=1$, which would contradict the monotony of $b_{i,\bullet}$. So, we can assume $j_1<j_2$ and $b_{i+1,j_1}-b_{i,j_1}=b_{i+1,j_2}-b_{i,j_2}=1$. Since the monotone sequences $b_{i,\bullet}$ and $b_{i+1,\bullet}$ both assume exactly three values, the only possibility that $j_2=j_1+1$, $b_{i+1,j_1}=b_{i+1,j_1+1}=w$ and $b_{i+1,j}=1$ for $j=i+2,i+3,...,j_1-1$ and $b_{i+1,j}=-1$ for $j=j_1+1,j_1+2,...,i$. Then we obtain $b_{i+1,j_1-1}-b_{i+1,j_1} = 1-w$, $b_{i,j_1-1}-b_{i,j_1} = 1-(w-1)=2-w$.
This implies that the smoothable edge $(j_1-1)\edge j_1$, has two lightly colored angles of the form $(j_1-1)\edge i\edge j_1$ and $(j_1-1)\edge (i+1)\edge j_1$. However, the equality $b_{i,j_1-1}-b_{i,j_1} = 1-w$ implies that $b_{j_1-1,j_1}=1-w$ and the equality $b_{i+1,j_1-1}-b_{i+1,j_1} = 2-w$ implies that $b_{j_1-1,j_1}=2-w$, a contradiction.

Under the assumptions of Subcase 2.1, so far we have proved that the sequence $b_{i+1,\bullet} $ assumes exactly three values $1$, $0$ and $-1$, and is non-increasing. Let us prove by induction on $s$ that each sequence $b_{s,\bullet}=\{b_{s,s+1}, b_{s,s+2}, ..., b_{s,s-1}\}$ assumes exactly three values $1$, $0$ and $-1$. Assuming this is true for $s=s_0$, let us prove it for $s=s_0+1$. Since $b_{s_0+1,j}- b_{s_0,j}\in \{0,1\}$ for $j\not=s_0,s_0+1$, we deduce that the sequence $b_{s_0+1,\bullet}$ assumes values from the set $\{-1,0,1,2\}$. Also, recall that $b_{s_0+1,\bullet}$ is non-increasing and assumes exactly three values by the assumption of Case 2. One of the three values $b_{s_0+1,\bullet}$ assumes has to be $-1$, otherwise the sum its entries will not equal zero. We claim that $b_{s_0+1,\bullet}$ has to contain a zero as well. To prove this, let us consider  the cases of even $\n$ and odd $\n$ separately.  If $\n$ is odd, $b_{s_0,\bullet}$ contains at least two zeros, otherwise the sum of its element could not be zero. Therefore, $b_{s_0+1,\bullet}$ contains a zero, unless $b_{s_0,\bullet}$ had exactly two zeros $b_{s_0, j_1}, b_{s_0,j_1+1}$ and both $b_{s_0+1,j_1}$, $b_{s_0,j_1+1}$ equal $1$. However, in that case the sequence $b_{s_0+1,\bullet}$ would have only $1$'s and $-1$'s, which would contradict the assumption of Case 2. If $\n$ is even, then $b_{s_0,\bullet}$ contains at least one zero, say $b_{s_0,j_1}=0$. Therefore, $b_{s_0+1,\bullet}$ contains a zero, unless one of the two colored angles of the edge $s_0\edge (s_0+1)$ is $s_0\edge j_1\edge (s_0+1)$. If the second colored angle is $s_0\edge (j_1+1)\edge (s_0+1)$, then $b_{s_0+1,j_1+1} = -1+1=0$; otherwise we obtain two contradictory values of the smoothable edge biresidue $b_{j_1,j_1+1}$.

This shows that two out of three values that $b_{s_0+1,\bullet}$ assumes are $-1$ and $0$.

Let us prove that $b_{s_0+1,\bullet}$ does not contain a two, that is, $b_{s_0+1,s_0+2}\not=2$. Assuming $b_{s_0+1,s_0+2}=2$, we obtain $b_{s_0,s_0+2}=1$ (the option $b_{s_0,s_0+2}=2$ is impossible, because $1=b_{s_0,s_0+1} \ge b_{s_0,s_0+2}$; the option $b_{s_0,s_0+2}=0$ is impossible, because it would lead to $b_{s_0,s_0+1}-b_{s_0,s_0+2} = \frac{1}{2} b_{s_0+1,s_0+2} $, a contradiction with smoothability of $(s_0+1)\edge (s_0+2)$). The biresidue $b_{s_0+1,s_0+3}$ equals either $0$, or $2$. The option $b_{s_0+1,s_0+3}=0$ would lead to $b_{s_0,s_0+3}=0$ and this would lead to two contradictory values of $b_{s_0+2,s_0+3}$; the option $b_{s_0+1,s_0+3}=2$ would lead to $b_{s_0,s_0+3}=1$, $b_{s_0+1,s_0+4}=0$, $b_{s_0,s_0+4}=0$, and this would lead to two contradictory values of $b_{s_0+3,s_0+4}$. 

Therefore, $b_{s_0+1,\bullet}$ assumes exactly three values $1$, $0$, $-1$. By induction, we argue that each $b_{s,\bullet}$ assumes exactly these three values. In particular, we get $b_{i,i+1}=1$ for each $i$. This reduces the proof to Case 1, which we have already discussed.

\textbf{Subcase 2.2.} Additionally to the assumption of Case 2, we assume that $b_{i,i+1}\not=b_{i+1,i+2}$ for each $i$.

Let $i$ be an index making $b_{i,i+1}$ the maximal possible. Then $b_{i,i+1}>b_{i+1,i+2}$. Without loss of generality, we may assume $b_{i,i+1}=2$. Then the angle $(i+1)\edge i\edge (i+2)$ must be lightly colored, and we obtain $b_{i+1,i+2}=1$. We claim that $\n$ must be even and that $b_{i+\n/2-1,i+\n/2}=2$.

Since the edge $(i+1)\edge (i+2)$ is smoothable, we get $b_{i,i+2}\le b_{i+1,i+2}=1$. If $b_{i,i+2}<1$, then smoothability of $i\edge (i+1)$ would imply that $b_{i,i+2}=-1$, which would make $\sum_j b_{i,j}=0$ impossible. Therefore, we have $b_{i,i+2}=1$. Therefore, two out of the three values the sequence $b_{i,\bullet}$ attains are $2$ and $1$. The third value must be of the form $1-\frac{2}{2^s}$, for some $s\ge 0$, because in the presence of lightly colored angles only \cite[Lemma 5.13]{Matviichuk2020} implies  that each $b_{i',i'+1}$ equals $\frac{2}{2^s}$ for some $s\ge 0$. If $s\ge 1$, then all the values of $b_{i,\bullet}$ will be non-negative, which would make $\sum_j b_{i,j}=0$ impossible. Therefore, $s=0$ and the three values $b_{i,\bullet} $ attains are $2$, $1$ and $-1$. Since there is precisely one $2$ in $b_{i,\bullet}$, the equality $\sum_j b_{i,j}=0$ implies that $\n$ is even and $b_{i,j}=1$ for  $i+2\le j \le i+ \n/2-1$ and $b_{i,j}=-1$ for $i+\n/2 \le j \le \n+i-1$. We obtain that $b_{i,i+\n/2-1}-b_{i,i+\n/2} = 2 \implies b_{i+\n/2-1,i+\n/2}=2$, as claimed.

Since $\n$ is even, we can either have $\n\equiv 0 {\rm ~mod~}4$ or $\n\equiv 2 {\rm~mod~}4$. In the former case, we have $\gcd(\n,\n/2-1)=1$. So, if $\n\equiv 0 {\rm ~mod~}4$, then by applying the claim above repeatedly, we obtain obtain $b_{i,i+1}=1$ for all $i$, a contradiction to the assumption of Subcase 2.2. Therefore, we must have $\n=2(2k+1)$ for some $k\in\mathbb{Z}_{\ge1}$. Then $\gcd(\n,\n/2-1)=2$, and applying the claim above repeatedly, we obtain $b_{i,i+1}=2$ for every other $i$. This implies that the biresidue matrix $B$ is equivalent to $\mathcal{Y}_{2(2k+1)}$.

\textbf{Case 3.}
There is a vertex $v$ such that there is only one colored angle of the form $i\edge v\edge (i+1)$ (possibly darkly colored).

The sequence $b_{i,\bullet}$ assumes only two values, say $\alpha>0$ and $-\beta<0$. Without loss of generality, we assume $\alpha \ge \beta$. We are going to show that $\alpha=\beta$, unless $\n=4$, in which case,  $B$ is projectively equivalent to $\mathcal{X}_4$. 

\textbf{Subcase 3.1.} $\alpha>\beta$.

\textit{Subcase 3.1.1.} $b_{i,i+1}=\alpha$, $b_{i,i+2}=-\beta$.

Smoothability of $i\edge (i+1)$ implies that $b_{i+1,i+2}$ is either $-\beta+2\alpha$ or $-\beta+\alpha$. 

If $b_{i+1,i+2}=-\beta+2\alpha$, then smoothability of $(i+1)\edge (i+2)$ implies that either $\alpha+\beta=2\alpha-\beta$ (in case the angle $(i+1)\edge i\edge (i+2)$ is lightly colored), or $\alpha+\beta=2(2\alpha-\beta)$ (in case the angle $(i+1)\edge i\edge (i+2)$ is darkly colored). If $\alpha+\beta=2(2\alpha-\beta)$, then $\alpha=\beta$, a contradiction. If $\alpha+\beta=2\alpha-\beta$, then $\alpha=2\beta$. Then we must have $\n=4$, the $i$-th row of the biresidue matrix, up to a multiple and a cyclic permutation, equals $(0,2,-1,-1)$. This leads to $B$ being projectively equivalent  to $\mathcal{X}_4$.

If $b_{i+1,i+2}=-\beta+\alpha$, then smoothability of $(i+1)\edge (i+2)$ implies that either $\beta+\alpha=\alpha-\beta$ (in case the angle $(i+1)\edge i\edge (i+2)$ is lightly colored), or $\beta+\alpha=2(\alpha-\beta)$ (in case the angle $(i+1)\edge i\edge (i+2)$ is darkly colored). If $\beta+\alpha=\alpha-\beta$ then $\beta=0$, a contradiction. If $\beta+\alpha=2(\alpha-\beta)$, then $\alpha=3\beta$. Then we must have $\n=5$, and the $i$-th row of the biresidue matrix, up to a multiple and a cyclic shift, equals $(0,3,-1,-1,-1)$. This quickly leads to a contradiction.

\textit{Subcase 3.1.2.} $b_{i,i+1}=\alpha$, $b_{i,i+2}=\alpha$, $b_{i,i+2}=-\beta$.

Smoothability of the edge $i\edge (i+1)$ implies that the value of $b_{i+1,i+2}$ is either $\alpha$, or $2\alpha$, or $3\alpha$. 

If $b_{i+1,i+2}=3\alpha$, then $b_{i+1,i+3}=-\beta$. By considering the four biresidues $b_{i,i+2}>b_{i,i+3},b_{i+1,i+2}>b_{i+1,i+3}$, we obtain two contradictory values of the smoothable edge biresidue $b_{i+2,i+3}$. 

If $b_{i+1,i+2}=2\alpha$, then $b_{i+1,i+3}=-\beta+\alpha$. Considering the lightly colored angle $(i+2)\edge i\edge (i+3)$, we obtain $b_{i+2,i+3}=\beta+\alpha$. Since we have assumed $\alpha>\beta$, we get $b_{i+1,i+3}-b_{i+2,i+3} = 2\beta < 2\alpha= b_{i+1,i+2}$, a contradiction to smoothability of the edge $(i+1)\edge (i+2)$.

If $b_{i+1,i+2}=\alpha$, then the angle $i\edge (i+3)\edge (i+1)$ has to be colored. Then $b_{i+1,i+3}$ equals either $-\beta+\alpha$ or $-\beta+2\alpha$. In either case, be must have $b_{i+1,i+2}=b_{i+1,i+3}$, otherwise the two colored angles $(i+2)\edge i\edge (i+3)$ and $(i+2)\edge (i+1)\edge (i+3)$ would lead to two contradictory values of $b_{i+2,i+3}$. Therefore, either $-\beta+\alpha = \alpha$, or $-\beta+2\alpha=\alpha$. In the former case we obtain $\beta=0$, a contradiction. The latter case leads to $\alpha=\beta$, a contradiction.

\textit{Subcase 3.1.3.} $b_{i,i+1}=b_{i,i+2}=b_{i,i+3}=\alpha$.

Let $j_1\ge3$ be the greatest $j$ such that $b_{i,i+j}=\alpha$. Then $b_{i,i+j_1+1}=-\beta$.

Smoothability of the edge $i\edge (i+1)$ implies that $b_{i+1,i+2}$ is either $\alpha$, or $2\alpha$ or $3\alpha$. 

If $b_{i+1,i+2}=3\alpha$, then the edge $i\edge (i+1)$ has only one darkly colored angle $i\edge (i+2)\edge (i+1)$. In particular, $b_{i+1,i+j_1}=b_{i,i+j_1}=\alpha$ and $b_{i+1,i+j_1+1}=b_{i,i+j_1+1}=-\beta$. Then the edge $(i+j_1)\edge (i+j_1+1)$ has two simply colored angles at $i$ and $i+1$. Therefore, we get $b_{i+2,i+j_1}=b_{i+2,i+j_1+1}$. Considering the smoothable edge $(i+1)\edge (i+2)$, we obtain that either $\alpha=-\beta + 3 \alpha$, or $\alpha = - \beta+6\alpha$. Both these options contradict $\alpha>\beta>0$.

If $b_{i+1,i+2}=2\alpha$, then the edge $i\edge (i+1)$ has a lightly colored angle $i\edge (i+2)\edge (i+1)$. The remaining lightly colored angle has to be either at $j_1+1$, or at $i+3$ (otherwise, it would contradict monotony of $b_{i+1,\bullet}$). If the remaining lightly colored angle for the edge $i\edge (i+1)$ is at $j_1+1$, then either $b_{i+1,j_1+1}=\alpha$, which would imply that $\alpha=\alpha-\beta$, a contradiction, or $-\beta <b_{i+1,j_1+1}<\alpha$, which would imply $\alpha+\beta = b_{j_1,j_1+1}=\alpha-b_{i+1,j_1+1}<\alpha+\beta$, a contradiction. If the remaining lightly colored angle for the edge $i\edge (i+1)$ is at $i+3$, then $j_1>3$ (otherwise we would get $\alpha+\beta = b_{i+3,i+4}= \alpha+2\beta$). Since the angles $j_1\edge i\edge (j_1+1)$ and $j_1\edge (i+1)\edge (j_1+1)$ are colored, the angle $j_1\edge (i+2)\edge (j_1+1)$ cannot be colored. Therefore, $b_{i+2,j_1}=\alpha=b_{i+2,j_1+1}$, and so either $\alpha= -\beta +2\alpha$, or $\alpha=-\beta+ 4\alpha$. Both these options contradict $\alpha>\beta>0$.

If $b_{i+1,i+2}=\alpha$, then monotony of $b_{i+1,\bullet}$ implies that the angle $i\edge (j_1+1)\edge (i+1)$ must be colored. So, either $-\beta<b_{i+1,j_1+1}<\alpha$, which would lead $\alpha+\beta = b_{j_1,j_1+1} = \alpha - b_{i+1,j_1+1}<\alpha+\beta$, a contradiction, or $b_{i+1,j_1+1}=\alpha$, which would lead to either $\alpha = -\beta+\alpha$ or $\alpha=-\beta+2\alpha$, both contradicting $\alpha>\beta>0$.

\textbf{Subcase 3.2.} $\alpha=\beta$.

By an overall rescaling, we can assume $\alpha=\beta=1$. Therefore the $i$-th row of the biresidue matrix is a cyclic shift of
$
(0,1,1,...,1,-1,-1,...,-1)
$.
The number of ones has to be equal to the number of negative ones, so $\n$ must be odd in Subcase 3.2.

The options for the value of $b_{i+1,i+2}$ are $1$, $2$ and $3$. The option $b_{i+1,i+2}=3$ leads to a contradiction, because monotony of $b_{i+2,\bullet}$ implies $b_{i+2,i+3}\in\{4,7\}$, both options  contradicting smoothability of the edge $(i+2)\edge (i+3)$. So, either $b_{i+1,i+2}=1$ or $b_{i+1,i+2}=2$. 

If $b_{i+1,i+2}=1$ then the $(i+1)$-th row must be a cyclic shift of the $i$-th row. 

If $b_{i+1,i+2}=2$, then the analysis of smoothability conditions shows that either $\n=5$ and the biresidue matrix is $\mathcal{X}_5$, or $\n>5$ and the rows $i,i+1,i+2,i+3,i+4,i+5$ have to obtained by a simultaneous cyclic shift from the following six rows:
$$
\begin{pmatrix}
0 & 1 & 1 & 1 & 1 & 1&1&...& 1 & -1 & -1& -1& -1&-1&-1&...& -1  \\
-1& 0 & 2 & 2 & 1 & 1&1&...& 1 & -1 & -1& -1& -1&-1&-1&...& -1 \\
-1&-2 & 0 & 2 & 1 & 1&1&...& 1 &  1 &  1& -1& -1&-1&-1&...&-1\\
-1&-2 &-2 & 0 & 1 & 1&1&...& 1 &  1 &  1&  1&  1&-1&-1&...&-1\\
-1 & -1 &-1&-1& 0 & 1&1&...& 1 &  1 &  1&  1&  1&-1&-1&...&-1\\
-1 & -1 & -1 &-1&-1& 0 & 1&...& 1 &  1 &  1&  1&  1& 1&-1&...&-1\\
\end{pmatrix}
$$

Therefore, the cyclic sequence $b_{i,i+1}, b_{i+1,i+2}, b_{i+2,i+3},...$ has groups of consecutive ones followed by groups of consecutive twos. Each group of ones consists of at least $3$ ones, and each group of twos consists precisely of $2$ twos. 

If there are no such groups of twos at all, that is $b_{l,l+1}=1$ for every $l$, then we easily obtain that $B=\mathcal{C}_{2k+1,k}=\mathcal{C}_{2k+1,k,\{0\}}$, where $\n=2k+1$. Therefore, we assume now that $b_{l,l+1}=2$ for some $l$.

Let us show the following periodicity condition:
 \begin{equation}\label{eq:periodicity_Zfamily}
 \left\{\begin{matrix}
 b_{l,l+1}=1, \\b_{l+1,l+2}=2,
 \end{matrix}\right. \xRightarrow{~~~~~~} 
  \left\{\begin{matrix}b_{l+(\n-5)/2,l+1+(\n-5)/2}=1, \\b_{l+1+(\n-5)/2,l+2+(\n-5)/2}=2.
   \end{matrix}\right. 
 \end{equation}
It is enough to show this for $l=i$.
Since $b_{i,i+(\n-1)/2}=1$, $b_{i,i+(\n-1)/2+1}=-1$ and $b_{i+1,i+(\n-1)/2}=1$, $b_{i+1,i+(\n-1)/2+1}=-1$ we deduce that  $b_{i+(\n-1)/2,i+(\n-1)/2+1}=2$. Similarly, since
$b_{i+3,i+(\n-1)/2+4}=1$, $b_{i+3,i+(\n-1)/2+5}=-1$ and $b_{i+4,i+(\n-1)/2+4}=1$, $b_{i+4,i+(\n-1)/2+1+5}=-1$ we deduce that the value of $b_{i+(\n-1)/2+4,i+(\n-1)/2+5}$ is also $2$.

The three consecutive values $b_{i+j,i+j+1}$, $j=(\n-1)/2+1, (\n-1)/2+2, (\n-1)/2 + 3$ are located between the $2$'s $b_{i+(\n-1)/2,i+(\n-1)/2+1}$ and $b_{i+(\n-1)/2+4,i+(\n-1)/2+5}$. So, these three values cannot be all $2$'s, which implies that they all have to equal $1$. Therefore, in particular, $b_{i+(\n-1)/2-1,i+(\n-1)/2}=2$, $b_{i+(\n-1)/2-2,i+(\n-1)/2-1}=1$, which proves \eqref{eq:periodicity_Zfamily}.

The periodicity condition \eqref{eq:periodicity_Zfamily} implies that $b_{i+1+sd}=2$ for each $s$, where $d=\gcd(\n,(\n-5)/2)$. We have $d=1$, unless $\n=5(2k+1)$ for some $k\in\mathbb{Z}_{\ge1}$. Since not all biresidues $b_{i,i+1}$ equal $2$, we obtain that $\n=5(2k+1)$ and $B$ is equivalent to $\mathcal{Z}_\n$.
\end{proof}

Recall that the rank of a skew-symmetric $\n\times \n$ matrix $B$ determines the smallest dimension $m$ such that $B$ can be realized as the biresidue matrix of a log symplectic form on $\mathbb{P}^{\n-1}\times \mathbb{D}^m$ with the toric polar divisor (see \autoref{lm:boundDimPolydisc}). Hence, it would be interesting to solve the following

\begin{problem}
Find the rank of each biresidue matrix in the list given in \autoref{thm:smoothableCycles}.
\end{problem}

It is easy to check that ${\rm rank}(\mathcal{X}_4) = 2$, ${\rm rank}(\mathcal{X}_5) = 2$, and
numerics seem to suggest that ${\rm rank}(\mathcal{Y}_\n) = \n-2$, ${\rm rank}(\mathcal{Z}_\n) = \n-3$ for admissible values of $\n$. However, the pattern for the rank of $\mathcal{C}_{\n,k,I}$ seems less clear. Nevertheless, we can prove the following partial results.

\begin{lemma} \label{lm:rankOfCNk}
The rank of $\mathcal{C}_{\n,k}$ is equal to the number of integers $0\le j \le \n-1$ such that neither $jk$, nor $j(k+1)$ is divisible by $\n$.
\begin{proof}
Denoting $\xi = e^{\frac{2\pi\ii}{\n}}$, we have the following eigenvectors for $\mathcal{C}_{\n,k}$:
$${\bf v}_j = \left(1, \xi^j, \xi^{2j}, \dots,  \xi^{(\n-1)j}\right), j=0,1,...,\n-1.$$
They are linearly independent by a Vandermonde determinant argument. To calculate the rank of $\mathcal{C}_{\n,k}$, we need to calculate how many of these eigenvectors have non-zero eigenvalues. The eigenvector ${\bf v}_0=(1,1,...,1)$ has zero eigenvalue for any $k$. For $j\not=0$, the eigenvector ${\bf v}_j$ has the eigenvalue
$$
 \sum_{\ell=1}^{k} \xi^{j\ell} - \sum_{\ell=\n-k}^{\n-1} \xi^{j\ell}=\sum_{\ell=1}^{k} \xi^{j\ell} - \sum_{\ell=1}^{k} \xi^{-j\ell} = \xi^j \frac{\xi^{jk}-1}{\xi^j-1} - \xi^{-j}\frac{\xi^{-jk}-1}{\xi^{-j}-1} = 
$$
$$=\xi^j \frac{\xi^{jk}-1}{\xi^j-1} - \xi^{-jk}\frac{1-\xi^{jk}}{1-\xi^{j}} = \frac{\xi^j (1-\xi^{-(j+1)k}) (1-\xi^{jk})}{1-\xi^j}.
$$
This eigenvalue equals zero if and only if either $jk$, or $j(k+1)$ is divisible by $\n$.
\end{proof}
\end{lemma}

The case when $B$ has corank one is especially important, since it leads to log symplectic forms on the compact space $\mathbb{P}^{\n-1}$ (see \autoref{lm:boundDimPolydisc}). Those can be classified as follows.

\begin{corollary}\label{cor:smoothableCyclesPn}
Let $B$ be a biresidue $\n\times \n$ matrix in which each row has zero sum. Let the rank of $B$ equal $\n-1$. Then $B$ has a full smoothable cycle if and only if it is projectively equivalent to $\mathcal{C}_{\n,k}$ for some $1\le k <\n/2$ such that $\gcd(\n,k)=\gcd(\n,k+1)=1$.
\end{corollary}

\begin{proof}
By \autoref{thm:smoothableCycles}, it is enough to check which of the listed there matrices have corank one. 

As we mentioned above ${\rm rank}(\mathcal{X}_4)=2$ and ${\rm rank}(\mathcal{X}_5)=2$, so these two have higher coranks. 

Each matrix $\mathcal{Y}_{2(2k+1)}$, $k\ge 1$, has even size. This together with skew-symmetry and the fact its rows add up to zero implies that the corank of $\mathcal{Y}_{2(2k+1)}$ is at least two.

Next, consider the case $B=\mathcal{Z}_{5(2k+1)}$, $k\ge 1$. Let ${\bf r}_0,\dots,{\bf r}_{\n-1}$ be the rows of $B$. Then one can check that, in addition to the relation $\sum_{j=0}^{\n-1}{\bf r}_j = {\bf 0}$, the rows satisfy the relation $\sum_{j=0}^{2k} \big({\bf r}_{5j} + {\bf r}_{5j+3} - {\bf r}_{5j+4}\big) = {\bf 0}$. This shows that the corank of $B$ in this case is greater than one.

Finally, let us consider the case $B=\mathcal{C}_{\n,k,I}$, where $1\le k<\n/2$, and $I=\{i_0,\dots,i_{d-1}\}$ is a sequence of zeros and ones of length $d=\gcd(\n,k)$. By replacing $k$ with $k-1$, if needed, we can assume that $i_s=0$ for some $0\le s <d$. By conjugating $B$ via a permutation matrix, we can assume $i_0=0$.
Let ${\bf r}_0,\dots, {\bf r}_{\n-1}$ be the rows of $B$, so that for instance ${\bf r}_0$ is ${\bf r}(k)$ as in \eqref{eq:row_basic}. By assumption, we have $\sum_{j=0}^{\n-1}{\bf r}_j={\bf 0}$. Additionally, if $d>1$, we have another relation $\sum_{j=0}^{\n/d-1} {\bf r}_{dj} = {\bf 0}$. Therefore, the corank of $B$ is greater than one, unless $d=1$, in which case we have $\mathcal{C}_{\n,k,I}=\mathcal{C}_{\n,k}$. Furthermore, if $\widetilde{d}=\gcd(\n,k+1)>1$, we have a new relation $\sum_{j=0}^{\n/\widetilde{d}-1} {\bf r}_{\widetilde{d}j} = {\bf 0}$. Therefore, the corank of $B$ is greater than one, unless $B$ is projectively equivalent to $\mathcal{C}_{\n,k}$ with $\gcd(\n,k)=\gcd(\n,k+1)=1$.  Conversely, \autoref{lm:rankOfCNk} implies that $ \mathcal{C}_{\n,k}$ has corank one provided that $\gcd(\n,k)=\gcd(\n,k+1)=1$.
\end{proof}

We finish this section with posing a problem about the following curious property of biresidue matrices. A biresidue matrix $B=(b_{ij})_{i,j=0}^{\n-1}$ is called \textit{holonomic}, if for every subset $J\subset\{0,1,\dots,\n-1\}$ of odd cardinality, the linear span of the rows of the submatrix $B_J = (b_{i,j})_{i,j\in J}$ does not contain the vector $(1,1,\dots,1)\in\mathbb{C}^J$. By \cite[Theorem 1.5]{Matviichuk2020}, holonomicity of $B$ is equivalent to holonomicity, in the sense of \cite{Pym2018}, of the corresponding log symplectic form $\omega_0$, and to finite-dimensionality of the Poisson cohomology of $\pi_0=\omega_0^{-1}$ everywhere locally. Numerical experiments seem to suggest the following

\begin{conjecture}
Every biresidue matrix listed in \autoref{thm:smoothableCycles} is holonomic.
\end{conjecture}

A result of similar nature has been proven in  \cite[Lemma 4.26]{Matviichuk2023} in the context of Hilbert schemes of Poisson surfaces. One should hope that the method of its proof can be adjusted to solve the current conjecture.

\newpage

\section{Examples of Poisson deformations produced by smoothable cycles}

\subsection{Feigin-Odesskii Poisson brackets}\label{subsec:FO}

Recall that for $1\le k<\n$ such that $\gcd(\n,k)=1$, one can define the Feigin-Odesskii bracket $q_{\n,k}$ on $\mathbb{C}^\n$ via the formula \cite[Section 5.2]{Hua2018}

\begin{equation}\label{eq:FO_explicit_formula}
    \{y_i,y_j\} = \left(
    \dfrac{\theta_{j-i}'(0)}{\theta_{j-i}(0)} + \dfrac{\theta'_{k(j-i)}(0)}{\theta_{k(j-i)}(0)} - 2\pi {\rm i} \n
    \right)
    y_iy_j + \sum_{r\not=0, j-i} \dfrac{\theta_{j-i+r(k-1)}(0) \theta_0'(0)}{\theta_{kr}(0)\theta_{j-i-r}(0)} y_{j-r} y_{i+r},
\end{equation}
where $i\not=j$ and the indices are considered modulo $\n$, and $\theta_\alpha$ denotes the theta function

$$
\theta_\alpha (\xi) = \theta(\xi+\frac{\alpha}{\n}\tau) \theta(\xi+\frac{1}{\n}+\frac{\alpha}{\n}\tau)~...~\theta(\xi+\frac{\n-1}{\n}+\frac{\alpha}{\n}\tau)  \exp\Big(2\pi {\rm i}(\alpha \xi +\frac{\alpha(\alpha-\n)}{2\n} \tau + \frac{\alpha}{2\n}) \Big),
$$
where
$$
\theta(\xi) = \sum_{k\in \mathbb{Z}} (-1)^k \exp\Big( 2\pi {\rm i} (k\xi + \frac{k(k-1)}{2} \tau)\Big),
$$
where $\tau \in \mathbb{C}$, ${\rm Im}(\tau)>0$, is a parameter defining an elliptic curve $E=\mathbb{C}/(\mathbb{Z}+\tau\mathbb{Z})$. Note that the Feigin-Odesskii bracket \eqref{eq:FO_explicit_formula} does not change when $\tau$ is replaced with $\tau+\n$. Therefore, it makes sense to parametrize the bracket \eqref{eq:FO_explicit_formula} by $\varepsilon = \exp(\frac{2\pi \ii}{\n} \tau)$. Let us consider its limit ${\rm Im}(\tau)\to +\infty$, or equivalently $\varepsilon\to0$. For an integer $\alpha$, we will use the notation $\overline{\alpha}$ for the remainder of $\alpha$ modulo $\n$. We will use the notation $k'$ for the unique integer such that $1\le k'<n$ and $kk' =1 {\rm~mod~} n$. The Poisson brackets $q_{n,k}$ and $q_{n,k'}$ are known to be isomorphic \cite[\S 3.4]{Chirvasitu2021}.

\begin{lemma}\label{lm:FO_toric}
The Feigin-Odesskii bracket \eqref{eq:FO_explicit_formula} has the following Taylor expansion when $\varepsilon \to 0$
\begin{equation}\label{eq:FO_asympt_expan}
q^0_{\n,k} + \varepsilon\, { C(\n)}\, q^{1}_{\n,k} + o \big( \varepsilon \big),
\end{equation}
where $C(\n)$ is a nonzero constant, and
$$
q^0_{\n,k} = \sum_{i=0}^{\n-1} \sum_{j=i+1}^{\n-1} 2\pi{\rm i} \Big({j-i} + \overline{k(j-i)} -\n\Big) y_i y_j ~\partial_{y_i} \wedge \partial_{y_j},
$$
$$
q^{1}_{\n,k} = \sum_{i=0}^{\n-1} y_{i+1}y_{i+k'}~ \partial_{y_i} \wedge \partial_{y_{i+k'+1}},
$$
where the indices are considered modulo $\n$.
\end{lemma}

\begin{proof}
Throughout the proof all the limits and Taylor expansions are meant as $\varepsilon\to 0$, or equivalently, as ${\rm Im}(\tau)\to +\infty$. Note that in the definition for $\theta(\xi)$ all the terms in the series tend to $0$, except the terms $k=0$ and $k=1$. Moreover, we obtain that
$$
\theta(\xi) = 1 - \exp(2\pi{\rm i} \xi) + o(\varepsilon).
$$
From here, we deduce that for $1\le \alpha<\n$, one has
$$
\theta_\alpha(\xi) = \exp\Big(2\pi {\rm i}(\alpha \xi +\frac{\alpha(\alpha-\n)}{2\n} \tau + \frac{\alpha}{2\n}) \Big) \prod_{j=0}^{\n-1} \left(1 - \exp\big(2\pi\ii(\xi+\frac{j}{\n})\big)\,\varepsilon^{\alpha} + o(\varepsilon)\right).
$$
Using that $\sum_{j=0}^{\n-1} \exp\big(2\pi\ii\frac{j}{\n}\big)=0$, we obtain
\begin{equation}\label{eq:thetaAlphaAsympt}
\theta_\alpha(\xi) =\exp\Big(2\pi {\rm i}(\alpha \xi +\frac{\alpha(\alpha-\n)}{2\n} \tau + \frac{\alpha}{2\n}) \Big) \big(1+ o(\varepsilon)\big), ~~~1\le \alpha<\n.
\end{equation}
Therefore, for any $1\le \alpha <\n$, we get
$$
\dfrac{\theta_\alpha'(\xi)}{\theta_\alpha(\xi)} =2\pi{\rm i}\alpha + o(\varepsilon).
$$

We deduce from this that for $0\le i<j\le \n-1$ one has

\begin{equation}\label{eq:1605136266}
\left(
    \dfrac{\theta_{j-i}'(0)}{\theta_{j-i}(0)} + \dfrac{\theta'_{k(j-i)}(0)}{\theta_{k(j-i)}(0)} - 2\pi {\rm i} \,\n
    \right) = 2\pi {\rm i} \Big( 
{j-i} + \overline{k(j-i)} -\n
\Big)
 + o(\varepsilon).
\end{equation}

Let us now show that the second summand on the right hand side of \eqref{eq:FO_explicit_formula} tends to $0$, and determine its leading term. To this end, let us define the nonzero constant $C(\n)$ via
$$ \theta'_0(0)=-2\pi {\rm i}\cdot \prod_{j=1}^{\n-1} \Big(1-\exp(2\pi{\rm i} \frac{j}{\n}) \Big) + o(\varepsilon) = C(\n) + o(\varepsilon). $$

Combining this with $\theta_0(0)=0$ and \eqref{eq:thetaAlphaAsympt}, we get for $r\not=0,j-i$ the following
\begin{equation}\label{eq:FOSecondTermAsympt}
\dfrac{\theta_{j-i+r(k-1)}(0)\theta'_0(0)}{\theta_{kr}(0)\theta_{j-i-r}(0)}  = \left\{
\begin{matrix}
0, & \text{ if } \overline{j-i+r(k-1)}=0, \\
\\
\pm C(\n)\varepsilon^{h\big(kr,j-i-r\big)}\big(1
 + o(\varepsilon)\big), & \text{otherwise}. 
\end{matrix} \right.
\end{equation}
where the function $h$ is defined by $h(\alpha,\beta) = g(\alpha)+g(\beta)-g(\alpha+\beta)$, $\alpha,\beta\in\mathbb{Z}$, where $g(\alpha) = \frac{1}{2}\overline{\alpha}(\n-\overline{\alpha})$, and the sign $\pm$ is $+$ if $\overline{kr} +\overline{j-i-r}<\n$, and $-$ otherwise.
 The following properties of the function $h(\alpha,\beta)$ are easily verified
\begin{enumerate}
\item\label{it:pos_h_N} $h(\alpha,\beta)$ is a positive integer for all integers $\alpha,\beta$ not divisible by $\n$.
\item $h(\alpha-1,\beta)<h(\alpha,\beta)$, if $\overline{\alpha}+\overline{\beta}\le \n$.
\item $h(\alpha+1,\beta)<h(\alpha,\beta)$, if $\overline{\alpha}+\overline{\beta}\ge \n$.
\item\label{it:min_h_N} $h(\alpha,\beta)=1$ if and only if either $\overline{\alpha}=\overline{\beta}=1$, or $\overline{\alpha}=\overline{\beta}=\n-1$.
\end{enumerate}

Property \ref{it:pos_h_N} shows that the second summand in the right hand side of \eqref{eq:FO_explicit_formula} converges to $0$, which together with \eqref{eq:1605136266} establishes the zeroth order term in the Taylor expansion \eqref{eq:FO_asympt_expan}. 

Property \ref{it:min_h_N} shows that the next term in the Taylor expansion \eqref{eq:FO_asympt_expan} is obtained from  the terms \eqref{eq:FOSecondTermAsympt} corresponding to either $\overline{kr}=\overline{j-i-r}=1$, or $\overline{kr}=\overline{j-i-r}=\n-1$. The former option establishes that the coefficient of $y_{i+1} y_{i+k'}\partial_{y_{i}}\wedge \partial_{y_{i+k'+1}} $ in the Taylor expansion of \eqref{eq:FO_explicit_formula} equals  $\varepsilon \, C(\n)+o(\varepsilon)$, and the latter option establishes the same about the coefficient of $y_{j+1} y_{j+k'}\partial_{y_{j}}\wedge \partial_{y_{j+k'+1}} $. 
\end{proof}

The Feigin-Odesskii Poisson bracket $q_{\n,k}$ on $\mathbb{C}^\n$ is quadratic, and therefore it defines a Poisson bracket on $\mathbb{P}^{\n-1}$, which we will denote again by $q_{\n,k}$. Results of \cite{Feigin1998} imply that $q_{\n,k}$ on $\mathbb{P}^{\n-1}$ is generically symplectic if and only if, in addition to $\gcd(\n,k)=1$, we have $\gcd(\n,k+1)=1$. The goal of this section is to show that each generically symplectic Poisson bracket $q_{\n,k}$ on $\mathbb{P}^{\n-1}$ can be obtained from a smoothing diagram $\mathcal{C}_{\n,\widetilde{k}}$ via  \autoref{thm:deforming_several_smootables}, for some $1\le \widetilde{k}<\n/2$.

\begin{proposition}\label{prop:FObijection}
Let $1\le k <\n$ be such that $\gcd(\n,k)=\gcd(\n,k+1)=1$. Then the toric Poisson bracket $q_{\n,k}^0$ on $\mathbb{P}^{\n-1}$ is log symplectic, and its biresidue matrix $B$ is projectively equivalent to $\mathcal{C}_{\n,\widetilde{k}}$, where
\begin{equation}\label{eq:bijection_explicit}
\widetilde{k}+1=\min \Big((1+k)^{-1}{\rm ~mod~}\n, (1+k^{-1})^{-1}{\rm ~mod~}\n \Big).
\end{equation}
Moreover, each summand in the expression for $q_{\n,k}^1$ equals $\rho_{ij}$ (as in \autoref{prop:smoothable_edge_deform}) for a smoothable edge $i\edge j$ of $B$.
\end{proposition}

\begin{proof}
Let ${\Pi}$ be the $\n\times \n$ coefficient matrix of the toric Poisson bracket $q_{\n,k}^0$, that is
$$
{\Pi}_{ij} = 2\pi\ii\Big(\overline{j-i} + \overline{k(j-i)} -\n\Big), ~~~~0\le i\not=j\le \n-1.
$$

Let $\widehat{\Pi}$ be the $(\n-1)\times (\n-1)$ coefficient matrix of the projectivization of $\frac{1}{2\pi\ii}q_{\n,k}^0$ in the affine chart $z_i=\frac{y_i}{y_0}$, $i=1,2,...,\n-1$. For $1\le i<j\le \n-1$, we have
$$
\widehat{\Pi}_{ij} = {\Pi}_{ij} - {\Pi}_{i0} - {\Pi}_{0j} = \left\{
\begin{matrix}
0, & \text{ if } \overline{k(j-i)}+\overline{ki}\ge \n, \\
-2\pi\ii \n, & \text{ if } \overline{k(j-i)}+\overline{ki}< \n.
\end{matrix}
\right.
$$

Let $B=-\frac{1}{2\pi \ii \n}P^{-1}\mathcal{C}_{\n,\widetilde{k}}P$, where $P=(P_{ij})_{i,j=0}^{\n-1}$ is the permutation matrix defined as 
$$
P_{ij} = \left\{
\begin{matrix}
1, & \text{ if } \overline{j} = \overline{(k'+1)i}, \\
0, & \text{ otherwise},
\end{matrix}
\right.
~~~~~~~\text{where}~k'=k^{-1}{\rm~mod~}\n.
$$

We claim that $B=(b_{ij})_{i,j=0}^{\n-1}$ is the biresidue matrix for $q_{\n,k}^0$. To prove this we need to show that $\widehat{\Pi}$ is the inverse of $\widehat{B}=(b_{ij})_{i,j=1}^{\n-1}$. Let us first calculate the product $B\Pi$. 

Let $S=(S_{ij})_{i,j=0}^{\n-1}$ be the permutation matrix such that $S_{ij} = 1$, if $\overline{j-i}=1$, and $S_{ij}=0$, otherwise. Note that $S^\n$ equals $I_\n$, the identity matrix of size $\n$. Additionally, we have $\Pi S=S\Pi$, $\mathcal{C}_{\n,\widetilde{k}} S=S\mathcal{C}_{\n,\widetilde{k}} $, $SP=P S^{k'+1}$, and $BS=SB$.

Calculating $B\Pi$ directly seems challenging; instead, we are going to calculate $(I_\n-S^{k'+1})B\Pi$. Tracing the definitions, we find that the $(i,j)$-entry of $(I_\n-S^{k'+1})B\Pi$ equals
$$
\dfrac{1}{2\pi \ii \n} \Big(
-\Pi_{ij}-\Pi_{i+k'+1,j}+\Pi_{i+1,j} + \Pi_{i+k',j}
\Big) = \left\{
\begin{matrix}
1, & \text{ if } i=j, \\
-1, & \text{ if } \overline{j-i}=k'+1, \\
0,& \text{ otherwise}.
\end{matrix}
\right.
$$
In other words, we obtain that $(I_\n-S^{k'+1})B\Pi = I_\n-S^{k'+1}$. Since the matrix $B\Pi- I_\n$ is annihilated after multiplication by $I_\n-S^{k'+1}$ on the right, each column of $B\Pi- I_\n$ is collinear with the constant vector $(1,1,...,1)$. Meanwhile, the matrix $B\Pi- I_\n$ commutes with $S^{k'+1}$, so, by a similar argument, each row of $B\Pi- I_\n$ is collinear with the constant vector, too. Hence, we have $B\Pi = I_\n + a U_\n$, where $a\in \mathbb{Q}$, and $U_\n$ is the $\n\times \n$ matrix whose each entry equals $1$. To determine $a$, we note that the sum of entries in each column of $B\Pi$ equals zero. We therefore obtain 
$$
B\Pi = I_\n - \frac{1}{\n} U_\n.
$$

This readily implies that $\widehat{B}\widehat{\Pi}=I_{\n-1}$, as desired. 
Proof of the claim about $q_{\n,k}^1$ is now straightforward and is left to the reader.
\end{proof}
 
Let us summarize the results of this section. 

\begin{theorem}\label{thm:FOviaSmoothableCycles}
Among the matrices listed in \autoref{thm:smoothableCycles}, the ones that can be realized as biresidue matrices of toric log symplectic forms on $\mathbb{P}^{\n-1}$ are precisely $\mathcal{C}_{\n,\widetilde{k}}$ with $\gcd(\n,\widetilde{k})=\gcd(\n,\widetilde{k}+1)=1$. 
Applying \autoref{thm:deforming_several_smootables} to the cycle of smoothable edges of $\mathcal{C}_{\n,\widetilde{k}}$, one obtains the log symplectic Feigin-Odesskii bracket $q_{n,k}$, where $k=((1+\widetilde{k})^{-1}-1) {\rm~mod~}\n$.
\end{theorem}

\begin{proof} The first part follows from \autoref{cor:smoothableCyclesPn} by taking into account \autoref{lm:boundDimPolydisc}. The second part follows from \autoref{lm:FO_toric} and \autoref{prop:FObijection}. 
\end{proof}

The Feigin-Odesskii brackets are commonly referred to as elliptic brackets due to their definition and geometry being tied to an elliptic curve. We conclude this subsection by providing a supporting argument that demonstrates their ellipticity in the sense of \autoref{def:elliptic}, provided they are log symplectic. This argument is not entirely rigorous, as it relies on two claims made in \cite{Feigin1998}. These claims, while widely believed to be correct, have not yet been mathematically proven to our best knowledge. The first claim is that the Poisson bracket $q_{n,k}$ on $\mathbb{P}^{n-1}$ given by \eqref{eq:FO_explicit_formula} has a modular-theoretic description by identifying $\mathbb{P}^{n-1}$ with $\mathbb{P}({\rm Ext}^1_E(\xi_{n,k},\mathcal{O}_E))$, where $\xi_{n,k}$ is an indecomposable rank $k$ degree $n$ bundle on $E$, unique up to fixing its determinant, and $\mathcal{O}_E$ is the trivial line bundle on $E$ (\cite[Section 5]{Hua2018} proves this claim for $k=1$). The second claim is that the symplectic leaves in this interpretation are given by the set of extensions 
\begin{equation}\label{eq:extensionOnE}
\xymatrix{
0\ar[r] & \mathcal{O}_E \ar[r]& V \ar[r] & \xi_{n,k}\ar[r] & 0},
\end{equation}
with a fixed isomorphism class of $V$ (\cite{Chirvasitu2023a} proves this for $k=1$). Let us accept both these claims.

The corank of the Poisson tensor at an extension \eqref{eq:extensionOnE} is equal to $\dim H^0(E, End(V))-1$ (\cite[Lemma 3.2]{Polishchuk1998}). From now on, we assume $\gcd(n,k)=\gcd(n,k+1)=1$. Using \cite{Atiyah1957}, it is easy to check that the unique open leaf is given by extensions with $V\cong \xi_{n,k+1}$. The codimension two leaves are parametrized by the rank $k+1$ degree $n$ vector bundles $V$ satisfying $\dim H^0(E, End(V)) = 3$. Such a $V$ has to have precisely two indecomposable summands $\xi_{n_1,k_1}$ and $\xi_{n_2,k_2}$ such that $n_1k_2 -n_2k_1 = \pm 1$. \autoref{lm:uniqueCodim2} below shows that these numbers $n_1,n_2,k_1,k_2$ are uniquely determined, up to swapping the indecomposable summands. We can choose the determinant of $\xi_{n_1,k_1}$ arbitrarily and the determinant of $\xi_{n_2,k_2}$ is then uniquely determined. Therefore, the codimension two symplectic leaves are parametrized by the elliptic curve $E$.

\begin{lemma}\label{lm:uniqueCodim2}
Let $1\le k <n$ be such that $\gcd(n,k+1)=1$. Then there are unique positive integers $n_1,n_2, k_1,k_2$ such that $n_1+n_2=n$, $k_1+k_2=k+1$  and $n_1k_2-n_2k_1 = 1$.
\end{lemma}

\begin{proof} \textit{Existence.} For a sequence of positive integers $q_1,q_2,...,q_g$, one can define the continued fraction as follows
$$
[q_1,q_2,...,q_g] := q_1 - \frac{1}{q_2-\frac{1}{\ddots q_{g-1}-\frac{1}{q_g}}}.
$$

One can express $\frac{n}{k+1}=[q_1,q_2,...,q_g]$ for some $q_1,q_2,...,q_g\ge 2$. Let us define $A\in {\rm SL}(2,\mathbb{Z})$ via
$$
A = T^{q_1} S T^{q_2} S \dots S T^{q_g},
\hskip0.1\textwidth
\text{where~}
S = \begin{pmatrix}
0 & -1 \\
1 & 0
\end{pmatrix},
~~
T = \begin{pmatrix}
1 & 1 \\
0 & 1
\end{pmatrix}.
$$
We have $A\begin{pmatrix}0\\1\end{pmatrix}=\begin{pmatrix}n\\k+1\end{pmatrix}$. Therefore, for $B = A T^{-1}$, we have $B\begin{pmatrix}1\\1\end{pmatrix}=\begin{pmatrix}n\\k+1\end{pmatrix}$. We now define the desired numbers $n_1,n_2,k_1,k_2$ using the entries of $B\in {\rm SL}(2,\mathbb{Z})$ as follows
$$
B = \begin{pmatrix}
n_1 & n_2 \\
k_1 & k_2
\end{pmatrix}.
$$
In terms of the continued fractions we have 
$$
\frac{n_1}{k_1} = [q_1,q_2,...,q_{g-1}],~~~~ \frac{n_2}{k_2} = [q_1,q_2,...,q_{g-1}, q_g-1].
$$

\textit{Uniqueness.} Let $\widehat{n}_1,\widehat{n}_2, \widehat{k}_1,\widehat{k}_2$ be another quadruple of positive integers satisfying the same conditions. Defining $B$ as above, and $\widehat{B}$ in a similar way using the hatted quadruple, we obtain that the matrix $B^{-1}\widehat{B}\in{\rm SL}(2,\mathbb{Z})$ satisfies $B^{-1}\widehat{B}\begin{pmatrix}1\\1\end{pmatrix} = \begin{pmatrix}1\\1\end{pmatrix}$. Such a matrix necessarily has to be of the form
$$
B^{-1} \widehat{B} = \begin{pmatrix} b +1 & -b \\ b & 1-b\end{pmatrix}, ~~\text{for some } b\in\mathbb{Z}.
$$
It follows that
$$
\begin{pmatrix} \widehat{n}_1&\widehat{n}_2 \\ \widehat{k}_1&\widehat{k}_2\end{pmatrix} = \widehat{B} = \begin{pmatrix}
n_1 & n_2 \\
k_1 & k_2
\end{pmatrix} + b \begin{pmatrix}
n_1+n_2 & -n_1-n_2  \\
k_1 + k_2 & -k_1 - k_2 
\end{pmatrix}.
$$
The conditions $k_1,k_2\ge1$ and $\widehat{k}_1,\widehat{k}_2\ge1$ now imply that $b=0$.
\end{proof}

\newpage

\subsection{Poisson deformation $\mathcal{C}_{4,1}$}\label{subsec:C41deform}

The biresidue matrix $\mathcal{C}_{4,1}$ and its smoothing diagram are as follows

$$
\begin{pmatrix}
 0 & 1 & 0 & -1 \\
-1 & 0 & 1 & 0 \\
 0 &-1 & 0 & 1 \\
 1 & 0 &-1 & 0
\end{pmatrix}
\hskip0.3\textwidth
\sageNgon[1.0]{4}{\smoothedge{v1}{v2}{v3}{v4} \smoothedge{v2}{v3}{v4}{v1}
\smoothedge{v3}{v4}{v1}{v2}
\smoothedge{v4}{v1}{v2}{v3}}
$$

Since the rank of $\mathcal{C}_{4,1}$ is two, we can realize this biresidue matrix on $\mathbb{P}^3\times \mathbb{D}^1$ (\autoref{lm:boundDimPolydisc}). Concretely, we can take the following log symplectic form
$$
\omega_0 = \frac{1}{4} \left(\frac{dz_1}{z_1}\wedge \frac{dz_2}{z_2} +\frac{dz_2}{z_2} \wedge \frac{dz_3}{z_3}\right) - \frac{1}{8} \left(\frac{dz_1}{z_1}-\frac{dz_2}{z_2} + \frac{dz_3}{z_3}\right) \wedge dx,
$$
where $z_i=\frac{y_i}{y_0}$, $i=1,2,3$, are affine coordinates on $\mathbb{P}^3$, and $x$ is the standard coordinate on $\mathbb{D}^1$. Let us invert $\omega_0$ to obtain
$$
\pi_0 = -2(z_1\partial_{z_1} \wedge z_2 \partial_{z_2} + z_1\partial_{z_1} \wedge z_3 \partial_{z_3} +z_2\partial_{z_2} \wedge z_3 \partial_{z_3}) + 4(z_1\partial_{z_1}+z_3\partial_{z_3})\wedge \partial_x.
$$
Using the homogeneous coordinates $y_i$, $i=0,1,2,3$, on $\mathbb{P}^3$ we get
$$
\pi_0 = -\sum_{i=0}^3 y_i\partial_{y_i}\wedge y_{i+1} \partial_{y_{i+1}} - 2\left( y_0\partial_{y_0}-y_1\partial_{y_1}+y_2\partial_{y_2}-y_3\partial_{y_3}\right)\wedge \partial_x.
$$
The log symplectic form $\omega_0$ has four smoothable edges $i\edge (i+1)\text{~mod~} 4$, $i=0,1,2,3$. The corresponding individual Poisson deformations are given by
$$
\begin{matrix}
     \rho_{01} &=& e^{-x} ~y_2 y_3~ \partial_{y_0} \wedge \partial_{y_1},  \\
     \rho_{12} &=& e^{x}~ y_3 y_0 ~\partial_{y_1} \wedge \partial_{y_2}, \\
     \rho_{23} &=& e^{-x}~  y_0 y_1 ~\partial_{y_2} \wedge \partial_{y_3}, \\
     \rho_{30} &=& e^{x}~ y_1 y_2 ~\partial_{y_3} \wedge \partial_{y_0}.
\end{matrix}
$$

Let us find a Poisson deformation $\pi(\varepsilon) = \pi_0 + \varepsilon \pi_1 + \varepsilon^2 \pi_2 + ... $ such that $\pi_1= \sum_{i=0}^3 \rho_{i,i+1}$. We are going to look for $\pi(\varepsilon)$ of the following form
$$
\pi(\varepsilon) = \pi_0 + \varepsilon\, \sum_{i=0}^3 \xi(x,\varepsilon)^{(-1)^i}  \rho_{i,i+1} + \eta(x,\varepsilon) \sum_{i=0}^3 y_i\partial_{y_i} \wedge y_{i+1} \partial_{y_{i+1}}.
$$
After working out the Schouten integrability $[\pi(\varepsilon),\pi(\varepsilon)]=0$, one finds that the deformation $\pi(\varepsilon)$ is Poisson if and only if $\xi,\eta$ satisfy the following system of ODEs (all derivatives are with respect to $x$):\begin{equation*}
\begin{aligned}
	\xi' &=& \xi \eta, ~~~~~~~~~~~~~~~~~~ &\\
	\eta'&=& \frac{\varepsilon^2}{4} \left(
	\frac{\xi^2}{e^{2x}} -\frac{e^{2x}}{\xi^{2}}
	\right).
\end{aligned} 
\end{equation*}
The Picard-Lindel\"{o}f theorem implies, that one can solve these ODEs uniquely, after specifying the initial conditions at $x=0$, and possibly shrinking the disc $\mathbb{D}^1$. For instance, we can specify $\xi\big|_{x=0} = 1$, $\eta\big|_{x=0} = 0$.

Let us now find where the Poisson structure $\pi=\pi(\varepsilon)$ on $\X=\mathbb{P}^3\times \mathbb{D}^1$ drops rank. We calculate its Pfaffian
$$
\displaystyle \pi \wedge \pi\wedge \sum_{i=0}^3 y_i\partial_{y_i} = {\rm Pf}(\pi) \bigwedge_{i=0}^3 \partial_{y_i} \wedge \partial_x,
$$
$$
{\rm Pf}(\pi) = 32(-1+\eta)y_0y_1y_2y_3 + 8\varepsilon \left(
y_2^2 y_3^2 \frac{\xi}{e^x} + 
y_3^2 y_0^2 \frac{e^x}{\xi} + 
y_0^2 y_1^2 \frac{\xi}{e^x} + 
y_1^2 y_2^2 \frac{e^x}{\xi} 
\right).
$$

Applying the Jacobian criterion to $\D = \{{\rm Pf}(\pi)=0\}$, we obtain that for small values of $\varepsilon\not=0$ and $x$, the divisor $\D$ is smooth away from $L_1 \times \mathbb{D}^1$ and $L_2\times \mathbb{D}^1$, where $L_1$ is the line $\{y_0=y_2=0\}\subset \mathbb{P}^3$ and $L_2$ is  the line $\{y_1=y_3=0\}\subset \mathbb{P}^3$. Similarly, the surface $\D \cap \{x=0\}$ is smooth away from $L_1$ and $L_2$. One can check that at each point on the lines $L_1$ and $L_2$, the surface $\mathsf{D}\cap \{x=0\}$ has a normal crossings singularity, except for four points on each line, where the singularity is the Whitney umbrella.

The space of codimension $2$ symplectic leaves is isomorphic to the quotient of the smooth surface $(\D\cap \{x=0\}) \setminus \{L_1,L_2\}$ by the orbits of the Hamiltonian vector field 
$$
u = [\pi, x] = [\pi_0,x] = -2 \Big( y_0 \partial_{y_0} - y_1\partial_{y_1} + y_2 \partial_{y_2} -y_3\partial_{y_3}\Big).
$$

\begin{figure}[h]
\centering
\includegraphics[scale=0.8]{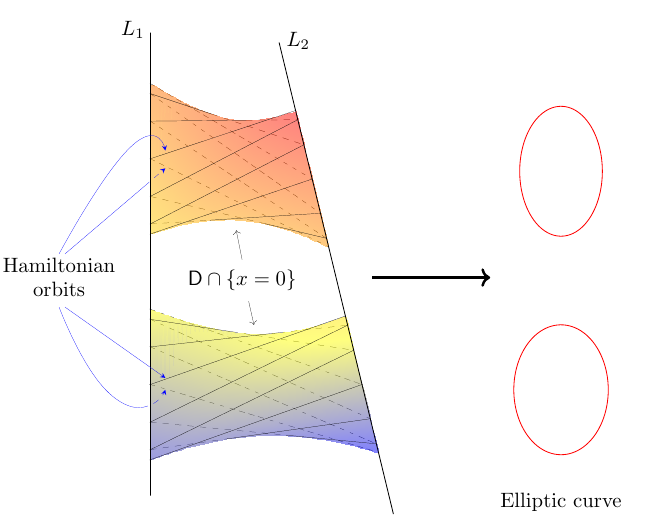}
\caption{Degeneracy divisor for the Poisson deformation $\mathcal{C}_{4,1}$.}
\label{fig:C41deform}
\end{figure}

Each orbit of $u$ is a straight line connecting a point of $L_1$ with a point of $L_2$. The quotient of $\mathbb{P}^3\setminus(L_1\cup L_2)$ by the Hamiltonian orbits is naturally isomorphic to $L_1\times L_2\cong \mathbb{P}^1\times \mathbb{P}^1$. The polynomial ${\rm Pf}(\pi)$ can be interpreted as a section of the anticanonical line bundle $\mathcal{O}_{\mathbb{P}^1\times\mathbb{P}^1}(2,2)$. Therefore, by the adjunction formula, the zero locus of $\mathbb{\rm Pf}(\pi)$ on $\mathbb{P}^1\times\mathbb{P}^1$ has trivial canonical bundle. This way we obtain that the space of codimension $2$ leaves of $\pi(\varepsilon)$ is an elliptic curve, for small non-zero values of $\varepsilon$. We illustrate this on \autoref{fig:C41deform}, where the surface $\D\cap\{x=0\}\subset\mathbb{P}^3$ takes on the shape of a pair of pillowcases.

\newpage

\subsection{Poisson deformation $\mathcal{X}_5$}\label{subsec:X5deform}

The biresidue matrix $\mathcal{X}_5$ and its smoothing diagram are as follows 
$$
(b_{ij})_{i,j=0}^4 = 
\begin{pmatrix}
 0 & 1 & 1 &-1 &-1 \\
-1 & 0 & 2 & 0 &-1 \\
-1 &-2 & 0 & 2 & 1 \\
 1 & 0 &-2 & 0 & 1 \\
 1 & 1 &-1 &-1 & 0
\end{pmatrix}
\hskip0.3\textwidth
\sageNgon[1.2]{5}{
\smoothedge{v1}{v2}{v3}{v4} 
\smoothedge{v2}{v3}{v4}{v5}
\smoothedge{v3}{v4}{v1}{v2}
\smoothedge{v4}{v5}{v2}{v3}
\smoothedge{v5}{v1}{v3}{v3}
}
$$

Since the rank of $\mathcal{X}_{5}$ is two, we can realize this biresidue matrix on $\mathbb{P}^4\times \mathbb{D}^2$ (\autoref{lm:boundDimPolydisc}). Concretely, we can take the following log symplectic form
$$
\omega_0 = \dfrac{1}{60}\left(\sum_{1\le i<j\le4} b_{ij}\dfrac{dz_i}{z_i} \wedge \dfrac{dz_{j}}{z_{j}}
\right) - \dfrac{1}{30}\left(-3\dfrac{dz_1}{z_1}+2\dfrac{dz_2}{z_2}-3\dfrac{dz_3}{z_3}+2\dfrac{dz_4}{z_4}\right)\wedge {dx_1}
+\dfrac{1}{30}\left(-\dfrac{dz_1}{z_1}+\dfrac{dz_3}{z_3}-2\dfrac{dz_4}{z_4}\right)\wedge {dx_2},
$$
where $z_i=\frac{y_i}{y_0}$, $i=1,2,3,4$, are affine coordinates on $\mathbb{P}^4$, and $x_1$, $x_2$ are the standard coordinates on $\mathbb{D}^2$. Let us invert $\omega_0$ to obtain
$$
\pi_0 = -4 (2~z_1\partial_{z_1} \wedge z_2 \partial_{z_2} + 2~z_1 \partial_{z_1} \wedge z_3 \partial_{z_3} + z_1 \partial_{z_1} \wedge z_4 \partial_{z_4} + 4~z_2 \partial_{z_2} \wedge z_3 \partial_{z_3} + 3~z_2 \partial_{z_2} \wedge z_4 \partial_{z_4} + z_3 \partial_{z_3} \wedge z_4 \partial_{z_4}) - 
$$
$$
- 5 (z_1 \partial_{z_1} + z_3 \partial_{z_3})\wedge \partial_{x_1} + 3(3~z_1\partial_{z_1} + 2~z_2 \partial_{z_2} + z_3 \partial_{z_3} + 4~z_4 \partial_{z_4})\wedge \partial_{x_2}.
$$
Using the homogeneous coordinates $y_i$, $i=0,1,2,3,4$, on $\mathbb{P}^4$ we get
$$
\pi_0 = -4 \left(\sum_{0\le i<j\le 4} b_{ij}~ y_i\partial_{y_i} \wedge y_j \partial_{y_j}\right) + u_1 \wedge \partial_{x_1} + u_2 \wedge \partial_{x_2},
$$
$$
u_1 = 2~y_0 \partial_{y_0} - 3~y_1 \partial_{y_1} + 2~ y_2 \partial_{y_2} - 3~y_3 \partial_{y_3} + 2~y_4\partial_{y_4},
$$
$$
u_2 = -3 \big(2~y_0 \partial_{y_0} - y_1 \partial_{y_1} + y_3 \partial_{y_3} - 2~y_4\partial_{y_4}\big).
$$

The log symplectic form $\omega_0$ has four smoothable edges $i\edge (i+1)\text{~mod~} 5$, $i=0,1,2,3,4$. The corresponding individual Poisson deformations are given by
$$
\begin{matrix}
     \rho_{01} &=&e^{10x_1-6x_2}~ y_2y_3~\partial_{y_0}\wedge \partial_{y_1}, \\
     \rho_{12} &=& e^{-5x_1+x_2}~ y_3y_4~\partial_{y_1}\wedge \partial_{y_2},\\
     \rho_{23} &=& e^{5x_1+x_2}~ y_0y_1~\partial_{y_2}\wedge \partial_{y_3}, \\
     \rho_{34} &=& e^{-10x_1-6x_2}~y_1y_2~\partial_{y_3}\wedge \partial_{y_4},\\
     \rho_{40} &=& e^{8x_2}~y_2^2~\partial_{y_4}\wedge \partial_{y_0}.\\
\end{matrix}
$$

Let us find a Poisson deformation $\pi(\varepsilon) = \pi_0 + \varepsilon \pi_1 + \varepsilon^2 \pi_2 + ... $ such that $\pi_1= \sum_{i=0}^4 \rho_{i,i+1}$. We are going to look for $\pi(\varepsilon)$ satisfying the following ansatz
$$
\pi (\varepsilon)= \pi_0 + \sum_{i=0}^4 \xi_i(x_1,x_2,\varepsilon) \rho_{i,i+1} + \eta\left(x_1,x_2,\varepsilon\right) \Big( y_1 \partial_{y_1} \wedge y_2 \partial_{y_2} + y_2 \partial_{y_2} \wedge y_3 \partial_{y_3} + y_3 \partial_{y_3} \wedge y_1 \partial_{y_1}\Big).
$$

The Schouten integrability $[\pi(\varepsilon),\pi(\varepsilon)]=0$ is equivalent to the following system of PDEs
$$
\begin{matrix}
(\xi_0)'_{x_1} & = & -\xi_0 \eta + \frac{1}{2} e^{-15x_1+15x_2} \xi_1 \xi_4, \\
(\xi_1)'_{x_1} & = & 0,\\
(\xi_2)'_{x_1} & = & 0,\\
(\xi_3)'_{x_1} & = & \xi_3 \eta - \frac{1}{2} e^{15x_1+15x_2} \xi_2 \xi_4 \\
(\xi_4)'_{x_1} & = & 0,\\
\eta'_{x_1} & = & \frac{1}{4}e^{-15x_1-5x_2} \xi_1 \xi_3 - \frac{1}{4} e^{15x_1-5x_2} \xi_0 \xi_2,
\end{matrix}
\hskip0.1\textwidth
\begin{matrix}
(\xi_0)'_{x_2} & = &\frac{1}{3}\xi_0 \eta - \frac{1}{2} e^{-15x_1+15x_2} \xi_1 \xi_4, \\
(\xi_1)'_{x_2} & = & 0,\\
(\xi_2)'_{x_2} & = & 0,\\
(\xi_3)'_{x_2} & = & \frac{1}{3}\xi_3 \eta - \frac{1}{2} e^{15x_1+15x_2} \xi_2 \xi_4 \\
(\xi_4)'_{x_2} & = & -\frac{2}{3}\xi_4\eta + \frac{1}{3} e^{-20x_2} \xi_0 \xi_3,\\
\eta'_{x_2} & = & \frac{1}{12}e^{-15x_1-5x_2} \xi_1 \xi_3 + \frac{1}{12} e^{15x_1-5x_2} \xi_0 \xi_2.
\end{matrix}
$$

It is straightforward to check that this system is compatible, in the sense that the two ways of calculating $\frac{\partial^2}{\partial x_1 \partial x_2}$ produce identical expressions. Therefore, the system can be uniquely solved, once we have specified the initial conditions. For instance, we can choose $\xi_i\big|_{x_1=x_2=0} =\varepsilon$, $0\le i \le 4$, $\eta\big|_{x_1=x_2=0} = 0$.

Let us now find where the Poisson structure $\pi=\pi(\varepsilon)$ on $\X=\mathbb{P}^4\times \mathbb{D}^2$ drops rank. We calculate its Pfaffian
$$
\displaystyle \pi \wedge \pi\wedge \pi\wedge \sum_{i=0}^4 y_i\partial_{y_i} = {\rm Pf}(\pi) \bigwedge_{i=0}^4 \partial_{y_i} \wedge \partial_{x_1} \wedge \partial_{x_2},
$$
$$
{\rm Pf}(\pi) = 720(-15+\eta) y_0y_1y_2y_3y_4 +
$$
$$+ 180\Big(e^{10x_1-6x_2}\xi_0 ~y_2^2y_3^2y_4 +
2e^{-5x_1+x_2} \xi_1~ y_0 y_3^2 y_4^2 +
2e^{5x_1+x_2} \xi_2~ y_0^2 y_1^2 y_4 +
e^{-10x_1-6x_2} \xi_3~ y_0 y_1^2 y_2^2 +
e^{8x_2} \xi_4~ y_1 y_2^3 y_3\Big).
$$

\begin{figure}[h]
\centering
\includegraphics[scale=1.2]{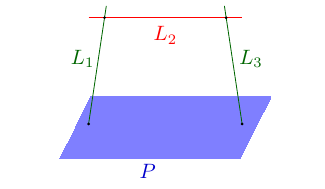}
%
%
\caption{Singular locus of the degeneracy divisor of the Poisson deformation $\mathcal{X}_5$}
\label{fig:X5Divisor}
\end{figure}

Applying the Jacobian criterion, we obtain that for small values of $x_1$, $x_2$ and $\varepsilon \not=0$, the degeneracy divisor $\mathsf{D}=\{{\rm Pf}(\pi)=0\} \subset \mathbb{P}^4 \times \mathbb{D}^2$ is smooth away from $L_j \times \mathbb{D}^2$, $j=1,2,3$, and $P\times \mathbb{D}^2$, where $L_j\subset \mathbb{P}^4$ are coordinate lines and $P\subset \mathbb{P}^4$ is a coordinate plane. Specifically, $L_1=\{y_0=y_2=y_3=0\}$, $L_2=\{y_0=y_2=y_4=0\}$, $L_3 = \{y_1=y_2=y_4=0\}$, and $P=\{y_1=y_3=0\}$. Likewise, the quintic threefold $\mathsf{Q}=\mathsf{D}\cap \{x_1=x_2=0\} \subset \mathbb{P}^4$ is smooth away from $L_1$, $L_2$, $L_3$, and $P$ (see \autoref{fig:X5Divisor}). The transverse slice of $\mathsf{Q}$ at a generic point of $L_2$ has an $\widetilde{E}_6$ singularity, while at a generic point of $L_1$ and $L_3$ it has an $\widetilde{E}_7$ singularity. The singularity of $\mathsf{Q}$ at a generic point of $P$ is normal crossings.

The space of codimension two symplectic leaves is isomorphic to the quotient of $\mathsf{Q}^\circ$, the smooth part of $\mathsf{Q}$, by the orbits of the Hamiltonian vector fields $u_1$ and $u_2$. Each orbit of $u_1$ is a straight line connecting a point of $L_2$ to a point of $P$. Hence, the quotient of $\mathbb{P}^4\setminus (L_2 \cup P)$ by $u_1$ is isomorphic to $L_2 \times P \cong \mathbb{P}^1 \times \mathbb{P}^2$. The intermediate quotient $\mathsf{Q}^\circ /u_1 \subset L_2\times P$ is the smooth part of an anticanonical divisor of $L_2\times P$, and its closure is singular at two points $S$ and $N$ obtained from $L_1$ and $L_3$, respectively. Therefore, the shape of this intermediate quotient is qualitatively similar to the ``onion'' depicted on \autoref{fig:X4Divisor}. The argument in the end of \autoref{ex:X4_full_deform} applies verbatim and shows that after taking the further quotient $\mathsf{Q}^\circ/u_1/u_2$ we obtain an elliptic curve.

\newpage
\bibliographystyle{hyperamsplain}

\end{document}